\documentclass{compositio}
\usepackage{enumitem}
\usepackage{amscd}
\usepackage{amsfonts}
\usepackage{graphicx}
\usepackage{amsmath}
\usepackage[all]{xy}
 \usepackage{amsmath}
\usepackage{color}
\usepackage{comment}
\usepackage{rlepsf}
 \usepackage{amssymb}

\newtheorem{Theorem}{Theorem}[section]
\newtheorem{Lemma}[Theorem]{Lemma}

\newtheorem{Definition}[Theorem]{Definition}
\newtheorem{Example}[Theorem]{Example}
\newtheorem{Proposition}[Theorem]{Proposition}
\newtheorem{Corollary}[Theorem]{Corollary}

\newtheorem{Remark}[Theorem]{Remark}

\numberwithin{equation}{section}

\newenvironment{Proof}[1][Proof]{\textbf{#1.} }{\ \rule{0.5em}{0.5em}}

\def\rref#1{(\ref{#1})}
\def \k {\kappa}
\def \Cl {\C_{\k}}
\def \ra {\xrightarrow}

\def \ad {\mathrm{ad}}

\def \Ad {\mathrm{Ad}}

\def \Aut {\rm Aut}

\def \C {\mathcal{C}}

\def \w {\omega}

\def \id {\mathrm{id}}
\def \ab{{\rm ab}}

\def \R {\mathbb{R}}

\def \X {\EuScript{X}}

\def \f {\phi}

\def \d {\partial}

\def \N {\mathbb N}

\def \tn {\otimes}

\def \tr {\triangleright}

\def \Gc {\mathcal{G}}

\def \ra {\xrightarrow}

\def \s {\scriptstyle}

\def \Z {\mathbb{Z}}

\def \tl {\triangleleft}

\def \p {\psi}

\def \fo {\textrm{ for each }}
\def \an {\textrm{ and }}

\def \wh {\textrm{ where }}
\def \de {\delta}

\def \GL {\rm GL}
\def \PGL {\rm PGL}

\def \Hom {\rm Hom}
\begin{document}

\title[Link invariants from finite categorical groups]{Link invariants from finite categorical groups and braided crossed modules}

\author{Jo\~{a}o Faria Martins}
\email{jn.martins@fct.unl.pt}
\address{  Departamento de Matem\'{a}tica and Centro de Matem\'{a}tica e Aplica\c{c}\~oes, Faculdade de Ci\^{e}ncias e Tecnologia, Universidade Nova de Lisboa, Quinta da Torre, 2829-516 Caparica, Portugal }
\author{Roger Picken}
\email{rpicken@math.ist.utl.pt}
\address{ Center for Mathematical Analysis, Geometry, and Dynamical Systems, {Departamento de Matem\'{a}tica, Instituto Superior T\'{e}cnico,} {  1049-001 Lisboa, Portugal }}

\thanks{ {J. Faria Martins was  partially supported by CMA/FCT/UNL, under the project PEst-OE/MAT/UI0297/2011. This work was partially supported by FCT (Portugal) through the projects
PTDC/MAT/098770/2008 
and PTDC/MAT/101503/2008. 
{We would like to thank Ronnie Brown for comments.}}}


\date{\today}



\begin{abstract}
 We define an invariant of tangles and framed tangles given a finite crossed module and a pair of functions, called a Reidemeister pair, satisfying natural properties. We give several examples of Reidemeister pairs derived from racks, quandles, rack and quandle cocycles, 2-crossed modules and braided crossed modules. We prove that our construction includes all rack and quandle cohomology (framed) link invariants, as well as the Eisermann invariant of knots, for which we also  find a lifting by using braided crossed  modules.
\end{abstract}

\maketitle

\section*{Introduction}
In knot theory, for a knot $K$, the fundamental group $\pi_1(C_K)$ of the knot complement $C_K$, also known as the knot group, is an important invariant, which however depends only on the homotopy type of the complement of $K$ (for which it is a complete invariant), and therefore, for example, it fails to distinguish between the square knot and the granny knot, which have homotopic, but not diffeomorphic complements. Nevertheless, a powerful knot invariant $I_G$ can be defined from any finite group $G$, by counting the number of morphisms from the knot group into $G$. In a recent advance, Eisermann \cite{E2} constructed, from any finite group $G$ and any $x \in G$, an invariant $E(K)$ that is closely associated to a complete invariant \cite{W}, known as the peripheral system, consisting of the knot group $\pi_1(C_K)$ and the homotopy classes of a meridian $m$ and a longitude $l$. Eisermann gives examples showing that his invariant is capable of distinguishing mutant knots  as well as detecting chiral (non-obversible), non-inversible and non-reversible  knots (using the terminology for knot symmetries employed in \cite{E2}). Explicitly the Eisermann invariant has the form, for a knot $K$:
$$E(K)=\sum_{\left\{\substack{ f\colon \pi_1(C_K) \to G\\ f(m)=x} \right\}} f(l),$$
and takes values in the group algebra $\Z[G]$ of $G$.

Eisermann's invariant has in common with many other invariants that it can be calculated by summing over all the different ways of colouring knot diagrams with algebraic data. Another well-known example of such an invariant is the invariant $I_G$ defined  above, which can be calculated by counting the number of colourings of the arcs of a diagram with elements of the  group $G$, subject to certain (Wirtinger) \cite{BZ} relations at each crossing.
 Another familiar construction is to use elements of a finite quandle to colour the arcs of a diagram, satisfying suitable rules at each crossing \cite{FR,CJKLS}. Note that the fundamental quandle of the knot complement is a powerful invariant that distinguishes all knots, up to simultaneous orientation reversal of $S^3$ and of the knot (knot inversion); see  \cite{J}. Other invariants refine the notion of colouring diagrams by assigning additional algebraic data to the crossings - a significant example is quandle cohomology \cite{CJKLS}. Eisermann's invariant can be viewed in several different ways, but for our purpose the most useful way is to see it as a quandle colouring invariant using a special quandle (the ``Eisermann quandle") associated topologically with the longitude and so-called partial longitudes coming from the diagram.

A diagram $D$ of a knot or link $K$ naturally gives rise to a particular presentation of the knot group, known as the Wirtinger presentation.
Our first observation is that there is also a natural crossed module of groups associated to a knot diagram \cite{BHS,BM,FM2}, namely $\Pi_2(X_D,Y_D)=\big(\partial\colon \pi_2(X_D,Y_D) \to \pi_1(Y_D)\big)$  - see the next section for the definition of a crossed module of groups and the description of $\Pi_2(X_D,Y_D)$.  This crossed module is a totally free crossed module \cite{BHS}, where  $\pi_1(Y_D)$ is the free group on the arcs of $D$ and   $\Pi_2(X_D,Y_D)$ is the free crossed module on the crossings of $D$.
The crossed module  $\Pi_2(X_D,Y_D)$ is not itself a knot  invariant, although it can be related to the knot group since $\pi_1(C_K)={\rm coker}(\d)$. However, up to crossed module homotopy \cite{BHS}, $\Pi_2(X_D,Y_D)$ is a knot invariant, depending only on the homotopy type of the complement $C_K$. Therefore, given a finite crossed module $\Gc=(\partial\colon E \to G)$, one can define a knot invariant $I_\Gc$  by counting all possible  colourings of the arcs and crossings of a diagram $D$ with elements of $G$ and $E$ respectively, satisfying some natural compatibility relations (so that colourings correspond to crossed module morphisms $\Pi_2(X_D,Y_D) \to \Gc$), and then normalising \cite{FM,FM2}. 

This invariant $I_\Gc(K)$  depends only on the homotopy type of the complement $C_K$ \cite{FM2,FM3}, thus it is a function of the knot group alone. Our main insight is that imposing a suitable restriction on the type of such colourings, and then counting the possibilities, gives a finer invariant. The restriction is to colour the arcs and crossings in a manner that is a) compatible with the crossed module structure, and b) such that the assignment at each crossing is given in terms of the assignments to two incoming arcs by two functions (one for each type of crossing), termed a Reidemeister pair. Since we are choosing particular free generators of $\Pi_2(X_D,Y_D)$ this takes away the homotopy invariance of the invariant. 

The two functions making up the Reidemeister pair must satisfy some conditions, and depending on the conditions imposed, our main theorem (Theorem \ref{maintheorem}) states that one obtains in this way an invariant either of knots or of framed knots (knotted ribbons). In fact our statement extends to tangles and framed tangles. 

This invariant turns out to have rich properties, which are described in the remainder of the paper (Section 4). It includes as special cases the invariants coming from rack and quandle colourings, from rack and quandle cohomology and the Eisermann invariant (subsections 4.1. and 4.2). Extending beyond these known examples, in subsection 4.3 we give a class of examples using 2-crossed modules \cite{Co}, which constitute a refinement of certain rack colouring invariants. In section 4.4. we introduce the notion of an Eisermann lifting, namely a Reidemeister pair derived from a braided crossed module \cite{BG} which reproduces the arc colourings of the Eisermann quandle, combined with additional information on the crossings. We give a simple example of an Eisermann lifting that is strictly stronger than the Eisermann invariant  it comes from. Finally, in subsection 4.5, we give a homotopy interpretation of the Eisermann liftings. For this we use the notion of non-abelian tensor product of groups, {defined by Brown and Loday {\cite{BrL0,BrL}.}}

\section{  Crossed modules of groups and the crossed module associated to a knot diagram}

\subsection{Definition of crossed modules and first examples}
\begin{Definition}[(Crossed module of groups)]\label{LCM}
 A crossed module of groups, {which we  denote by} ${\Gc= ( \d\colon E \to  G,\tr)}$, is given by a group morphism $\d\colon E \to G$ together with a  left action $\tr$ of $G$ on $E$ by automorphisms, {such that} the following conditions (called Peiffer equations) hold:
\begin{enumerate}
 \item $\d(g \tr e)=g \d(e)g^{-1}; \forall g \in G, \forall e \in E,$
  \item $\d(e) \tr f=efe^{-1};\forall e,f  \in E.$
\end{enumerate}
The crossed module of groups is said to be finite, if both groups $G$ and $E$ are finite. 
Morphisms of crossed modules are defined in the obvious way.
\end{Definition}

\begin{Example}
Any pair of finite groups $G$ and $E$, with $E$ abelian, gives a finite crossed module of groups with 
trivial $\d, \tr$ {(i.e. $\d(E)=1, \, g\tr e = e, \forall g \in G, \forall e \in E$). More generally we can choose any action of $G$ on $E$ by automorphisms, with trivial boundary map $\d\colon E \to G$.}
\label{simplestexample} 
\end{Example}

\begin{Example}
{Let $G$ be any finite group. Let $\Ad$ denote the adjoint action of $G$ on $G$. Then ${(\id \colon G \to G, \Ad)}$ is a finite crossed module of groups.}
\end{Example}

\begin{Example}\label{q}
Let $E$ be any finite group. Let $\Aut(E)$ be the group of automorphisms of $E$. Then $\Aut(E)$ is a finite group also.  Consider the map $\Ad\colon E \to \Aut(E)$ that sends $ e\in E$ to the automorphism $\Ad(e) \colon E \to E$. The group $\Aut(E)$ acts on $E$ as $\f \tr e=\f(e)$ where $\f\in \Aut(E)$ and $e \in E$. Then $({\Ad}\colon E \to \Aut(E),\tr)$ is a finite crossed module of groups.
\end{Example}

 \begin{Example}
There is a well-known construction of a crossed module of groups in algebraic topology, namely the fundamental crossed module associated to a pointed pair {$(X,Y)$ of path-connected topological spaces $(X,Y)$, thus $Y\subset X$, namely the crossed module:}
$$
\Pi_2(X,Y)= ( \d\colon \pi_2(X,Y) \to \pi_1(Y),\tr),
$$
with the obvious boundary map $\partial\colon \pi_2(X,Y) \to \pi_1(Y)$, {and the usual} action of $\pi_1(Y)$ on $\pi_2(X,Y)$; see figure \ref{action}, and  \cite{BHS} for a complete definition. This is an old result of Whitehead {\cite{W1,W2,W3}. }
\begin{figure} 
\centerline{\relabelbox 
\epsfysize 2.5cm
\epsfbox{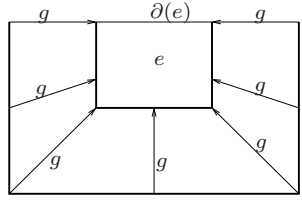}
\relabel{X}{$\s{g}$}
\relabel{Y}{$\s{g}$}
\relabel{Z}{$\s{g}$}
\relabel{W}{$\s{g}$}
\relabel{S}{$\s{g}$}
\relabel{T}{$\s{g}$}
\relabel{U}{$\s{g}$}
\relabel{e}{$\s{e}$}
\relabel{f}{$\s{\d(e)}$}
\endrelabelbox}
\caption{{\label{action} The action of an element $g \in \pi_1(Y)$ on an $e \in \pi_2(X,Y)$}.}
\end{figure}
\end{Example}

\begin{Example}
We may construct a topological pair $(X_D,Y_D)$ from a link diagram $D$ of a link $K$ {in $S^3$}, and thus obtain a crossed module $\Pi_2(X_D,Y_D)$, associated to the diagram. Regard the diagram as the orthogonal projection onto the $z=0$ plane in $S^3= \R^3 \cup \{\infty\}$ of a link $K_D$, isotopic to $K$, lying entirely in the plane $z=1$, except in the vicinity of each crossing point, where the undercrossing part of the link descends to height $z=-1$.
Then we take $X_D$ {(an excised link complement)} to be the link complement {$C_K$ of $K_D$ minus an open ball}, and $Y_D$ to be the $z\geq 0$ subset of $X_D$, i.e.
$$
X_D := \big (S^3 \setminus {n(K_D)}){\cap \{ (x,y,z)| z\geq -2\}}, \quad \quad Y_D := X_D \cap \{ (x,y,z)| z\geq 0\},
$$
{where $n(K_D)$ is an open regular neighbourhood of $K_D$ in $S^3$.} {Note that the space $X_D$ depends only on $K$ itself, so we can write it as $X_K$. The same is not true about $Y_D$.}

Each arc of the diagram $D$ corresponds to a generator of $\pi_1(Y_D)$ and there are no relations {between these}. Each crossing of the diagram $D$ corresponds to a generator of $\pi_2(X_D, Y_D)$, namely $\epsilon: [0,1]^2\rightarrow X_D$, where the image under $\epsilon$ of the interior of $[0,1]^2$ lies entirely in the region $z<0$ and the image of the boundary of $[0,1]^2$, a loop contained in $Y_D$, encircles the crossing as in Figure \ref{epsPQRS}.
\begin{figure}
\begin{center}
\includegraphics[width=10cm]{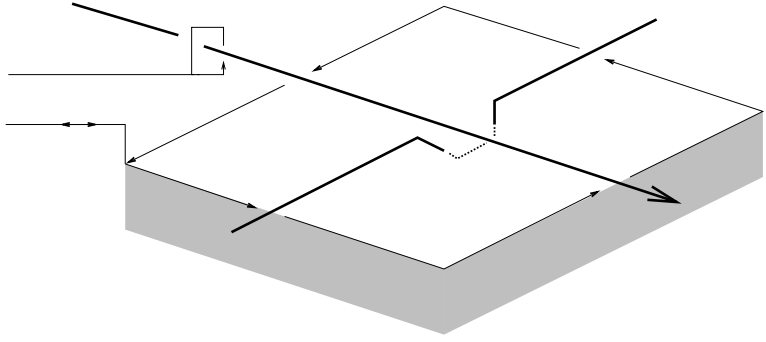}
\end{center}
\caption{\label{epsPQRS} {A generator of $\pi_1(Y_D)$ and a generator of $\pi_2(X_D, Y_D)$.}}
\end{figure}
This boundary loop is the product of four arc loops in $\pi_1(Y_D)$. 
{Again there are no {(crossed module)} relations between the generators of $\pi_2(X_D, Y_D)$ associated to the crossings}. (This  can  be justified by Whitehead's theorem \cite{W1,W2,W3, BHS}: for {path-connected spaces} $X, \,Y$, if $X$ is obtained from $Y$ by attaching 2-cells, then $\Pi_2(X,Y)$ is the free crossed module on the attaching maps of the 2-cells. Note that $X_D$ is homotopy equivalent to the CW-complex obtained from $Y_D$ by attaching a 2-cell for each crossing).

We observe that the quotient 
$\pi_1(Y_D) / {\rm im}\, \d$ is isomorphic to the fundamental group of the link complement {$C_K=\pi_1(S^3 \setminus n(K))$} for any diagram $D$, since quotienting $\pi_1(Y_D)$ by ${\rm im}\, \d$ corresponds to imposing the Wirtinger relations \cite{BZ}, which produces the Wirtinger presentation of {$\pi_1(C_K)$}, coming from the particular choice of diagram. Thus 
$\Pi_2(X_D,Y_D)$, whilst not being itself a link invariant {(unless considered up to crossed module homotopy)}, contains an important link invariant, namely {$\pi_1(C_K)$,} by taking the above quotient. The guiding principle in the construction to follow is to extract additional Reidemeister invariant information from the crossed module $\Pi_2(X_D,Y_D)$.
\label{fundexp}
\end{Example}

In the next subsection we will give a different perspective on the preceding example, via handle decompositions of the link complement in $S^3$ coming from a link diagram $D$, which allows us relate it to previous constructions of link invariants. This subsection is not necessary for following our construction in the rest of the paper, starting in section \ref{Xmod}.

\subsection{{Knot invariants from handle decompositions of (excised) link complements} }\label{elc}

Let $K$ be an oriented link in $S^3$.  Let $n(K)$ be an open regular neighbourhood of $K$ in $S^3$. {As above}, the link complement $C_K$ is defined as:
$$C_K=S^3 \setminus n(K) .$$

\begin{figure}
\begin{center}
\includegraphics[width=10cm]{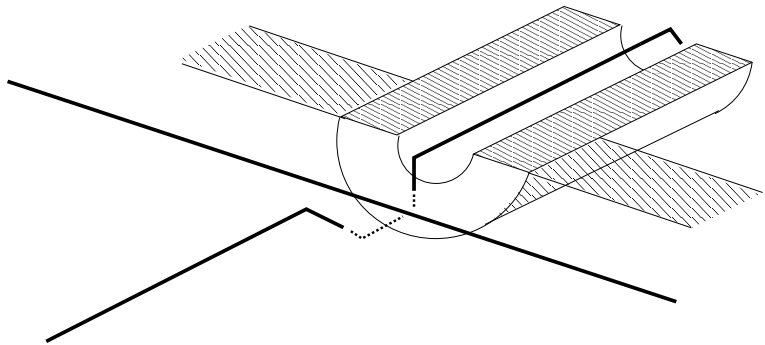}
\caption{\label{attach} The 1-handle associated with an arc of the projection.}
\end{center}
\end{figure}

\begin{figure}
\begin{center}
\includegraphics[width=7cm]{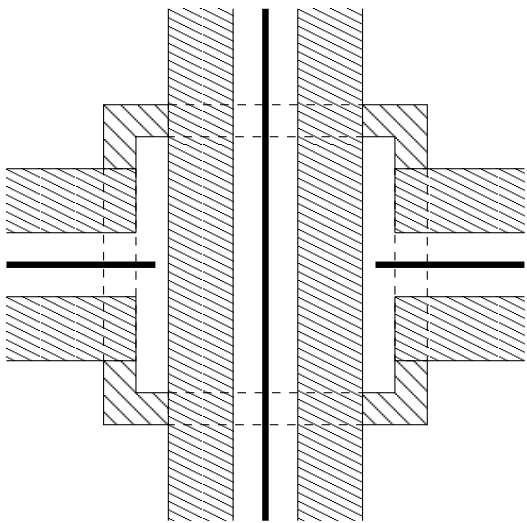}
\caption{\label{2handle} The 2-handle associated with a crossing of the projection.}
\end{center}
\end{figure}

Let us be given a regular {projection} of $K$, defining a diagram $D$ of $K$.
We can define a handle decomposition \cite{GS} of $C_K$ with a single $0$-handle. We have a 1-handle for each arc of the diagram $D$ of $K$, whose attaching region $I^2\times \d I$ is the more densely shaded region in figure \ref{attach}. Therefore each arc $c$ of the projection gives rise to an element $g_c \in \pi_1(C_K)$; these are the well known Wirtinger generators of $\pi_1(C_K)$. At each crossing $x$ we attach a 2-handle, whose attaching region $I\times \d I^2$ is the less densely shaded region shown (in part) in figure \ref{attach} and in figure \ref{2handle}. From this we immediately get the Wirtinger presentation of the fundamental group of a knot complement. {Note that attaching an $n$-handle is, up to homotopy, the same as attaching an $n$-cell}

Clearly, the union of the underlying 0-handle, the 1-handles associated to arcs of the projection and the 2-handles associated to crossings gives a handle decomposition of the space ${X_K=C_K \setminus {\rm int}( D^3)}$, where ${\rm int}( D^3)$ is the interior of the 3-ball, included in $S^3$ below the link $K$. To obtain $C_K$ from the {excised link} complement {$X_K$} we still have therefore to attach a 3-handle along the element {$t_D$ of {$\pi_2(X_K)$,}} associated with the boundary of $D^3$. 

Even though the handle decomposition of {${C}_K$ and of $X_K$} depends on the chosen diagram $D$ for $K$ it is proved in \cite{FM,FM2} that:
\begin{Theorem}
 Let $\Gc=(\d\colon E \to G, \tr)$ be a finite crossed module. Given a link $K$ with diagram $D$, let $\Pi_2(K,D)$ denote the fundamental crossed module of the pair $(C_K,C_K^1)$, where  $C_K^1$ is the handlebody made from the 0 and 1-handles of $C_K$. Then we have a knot invariant $I_\Gc$ where
$${I_\Gc(K)=\frac{1}{(\# G)^{\#\{\textrm{arcs of } D\}}} \# \Hom(\Pi_2(K,D), \Gc)\in \mathbb{Q}.}$$
{Here $\Hom(\Pi_2(K,D), \Gc)$ is the set of crossed module morphisms $\Pi_2(K,D) \to  \Gc$.}
\end{Theorem}
{As one might expect, this invariant $I_\Gc$ depends only on the homotopy type of $C_K$, and therefore only on $\pi_1(C_K)$. }

{Consider   the fundamental crossed module $\Pi_2(X_D,Y_D)$ of the pair {$( X_D,C_D^1)=( X_D,Y_D)$.} A modification of the invariant {${I}_\Gc$} above is to consider the finer knot invariant:}
$${\overline{I}_\Gc(K)=\frac{1}{(\# G)^{\#\{\textrm{arcs of } D\}}} \sum_{f \in \Hom({\Pi}_2(X_D,Y_D), \Gc) } f(t_D) \in \Z[\ker(\d\colon E \to G)].}$$

Let us give a nice geometric interpretation of $I_\Gc$ and $\bar{I}_\Gc$, coming from \cite{FM,FM2,FM3}. Let $\Gc=(\d\colon E \to G, \tr)$ be a finite crossed module. We can define its classifying space $B_\Gc$, thus $\pi_2(B_\Gc)=\ker(\d)$ {and  $\pi_1(B_\Gc)={\rm coker }(\d).$} Given a CW-complex $M$ let ${\rm TOP}(M,B_\Gc)$ be the topological space of maps $M \to B_\Gc$, where the topology is the $k$-iffication of the compact open topology in the space of continuous maps $M \to B_\Gc$. Therefore the set of homotopy classes of maps $M \to B_\Gc$ is exactly $\pi_0({\rm TOP}(M,B_\Gc))$.
\begin{Theorem}
 Given a link $K$:
$$I_\Gc(K)=\sum_{f \in \pi_0({\rm TOP}( C_K,B_\Gc))} \frac{1}{\# \pi_1({\rm TOP}( C_K,B_\Gc),f) } \in \mathbb{Q} $$
and
$$\bar{I}_\Gc(K)=\sum_{f \in \pi_0({\rm TOP}( {X_K,B_\Gc))}} \frac{1}{\# \pi_1({\rm TOP}( {X_K},B_\Gc),f) } {f(t_D)} \in \Z[\ker(\partial)].$$
\end{Theorem}
\begin{Proof}
Follows from the general construction in \cite{FM,FM2,FM3}.
\end{Proof}

\section{{Crossed modules of groups  and categorical groups }}\label{Xmod}

\subsection{A monoidal category $\C(\Gc)$ defined from a categorical group $\Gc$}\label{amc}
It is well-known that a crossed module of groups $\Gc$ gives rise to a categorical group, denoted $\mathcal{C}(\Gc)$, a monoidal groupoid where all objects and arrows are invertible, with respect to the tensor product operation; see \cite{BM,BHS,BL,FM,ML}. We recall the essential details.  Given a crossed module of groups ${\Gc= ( \d\colon E \to  G,\tr)}$, the {monoidal} category $\mathcal{C}(\Gc)$ has $G$ as its set of objects, and the morphisms from $U\in G$ to $V\in G$ are given by all pairs $(U,e)$ with $e\in E$ such that $\d(e)= VU^{-1}$. It is convenient to think of these morphisms as downward pointing arrows 
and / or  to represent them as squares - see \eqref{CGmorphisms}. 
\begin{equation}\label{CGmorphisms}
\xymatrix{
 & U \ar[d]_{(U,e) }\\ &V} \quad \quad \xymatrix{ & \ar@{-}[r]|U\ar@{-}[d] \ar@{-}@/_1pc/ @{{}{ }{}} [r]_{e}  
 & \ar@{-}[d] \\
                                                            & \ar@{-}[r]|V &} \begin{CD}\\\\ \textrm{ \quad \quad \quad with } \partial(e) U=V.\end{CD}
\end{equation}

The composition of morphisms $U\stackrel{e}{\rightarrow}V$ and 
$V\stackrel{f}{\rightarrow}W$ is defined to be $U\stackrel{fe}{\longrightarrow}W$, and the monoidal structure $\otimes$ is expressed as $U \tn V=UV$ on objects, and, on morphisms, as:
$$\begin{CD}
U\\@V(U,e)VV \\ V\end{CD} \,\,\otimes\,\, \quad  \quad \begin{CD} W\\@V{(W,f)}VV\\X\end{CD}   \quad = \, \begin{CD} UW \\ @VV{\big(UW,(V\tr f)\, e\big)}V\\ WX\end{CD}\,\,\,.
$$ 
These algebraic operations are shown using squares in  \eqref{ver} and \eqref{hor}. 
\begin{equation}\label{ver}
\xymatrix{& \ar@{-}[r]|U\ar@{-}[d] \ar@{-}@/_1pc/ @{{}{ }{}} [r]_{e} & \ar@{-}[d] \\& \ar@{-}[r]|V\ar@{-}[d] \ar@{-}@/_1pc/ @{{}{ }{}} [r]_{f} & \ar@{-}[d]\\& \ar@{-}[r]|W &
                                                            }\quad\quad \quad \xymatrix{ \quad \\ =}  {\xymatrix{ & \ar@{-}[r]|U\ar@{-}[dd] \ar@{-}@/_2pc/ @{{}{ }{}} [r]|{fe}  
 & \ar@{-}[dd] \\\\
                                                            & \ar@{-}[r]|W&}}
\end{equation}
\begin{equation}\label{hor}
 \xymatrix{ & \ar@{-}[r]|U\ar@{-}[d] \ar@{-}@/_1pc/ @{{}{ }{}} [r]_{e}  
 & \ar@{-}[d] \\
                                                            & \ar@{-}[r]|V &}
\quad \quad \quad {\begin{CD} \\ \\  \tn \end{CD} }  \xymatrix{ & \ar@{-}[r]|{U'}\ar@{-}[d] \ar@{-}@/_1pc/ @{{}{ }{}} [r]_{e'}  
 & \ar@{-}[d] \\
                                                            & \ar@{-}[r]|{V'} &} \quad \quad \quad { \begin{CD} \\ \\  = \end{CD} } 
 \xymatrix{ & \ar@{-}[rr]|{UU'}\ar@{-}[d] \ar@{-}@/_1pc/ @{{}{ }{}} [rr]_{{(V \tr e')\, e}}  &
 & \ar@{-}[d] \\
                                                            & \ar@{-}[rr]|{VV'} & &}
\end{equation}
The (strict) associativity of the composition and of the tensor product are trivial to check. The functoriality of the tensor product (also known as the interchange law) follows from the  second Peiffer equation. This calculation is done for example in {\cite{BM,BHS,Po}. }

Let $\k$ be any commutative ring. The monoidal category $\C(\Gc)$ has a $\k$-linear version $\Cl(\Gc)$, whose objects are the same as the objects of $\C(\Gc)$, but such that the set of morphisms $U \to V$ in $\Cl(\Gc)$ is given by the set of  all $\k$ linear combinations of morphisms $U \to V$ in $\C(\Gc)$. The composition and tensor product of morphisms in $\Cl(\Gc)$ are the obvious linear extensions of the  ones in $\C(\Gc)$. It is easy to see that $\Cl(\Gc)$ is a monoidal category.

The categorical group formalism is very well matched to the category of tangles to be used in the next section.

\section{Reidemeister $\Gc$ colourings of oriented tangle diagrams}
\subsection{Categories of tangles}
Tangles are a simultaneous generalization of braids and links. We follow \cite{O,K} very closely, to which we refer for more detais. Recall that an embedding of a manifold $T$ in a manifold $M$ is said to be neat if $\partial( T)=T \cap \partial( M)$.
\begin{Definition}
 An oriented tangle \cite{K,Tu,O} is a 1-dimensional smooth {oriented} manifold neatly embedded in $ \R \times \R\times [-1,1]$, such that $\partial(T) \subset \mathbb{N} \times\{0\} \times \{\pm 1\}$. A framed  oriented tangle is a tangle together with a choice of a framing in each of its components \cite{K,O}. Alternatively we can see a framed tangle as an embedding of a ribbon into $\R \times \R \times [-1,1]$, \cite{RT,Tu,CP}.
\end{Definition}
\begin{Definition}
 Two oriented tangles (framed oriented tangles) are said to be equivalent if they are related by an isotopy of $\R \times \R \times [0,1]$, relative to the boundary.
\end{Definition}

 \begin{Definition}A tangle diagram is a diagram of a tangle in $\R \times [-1,1]$, obtained from a tangle by projecting it onto $\R \times \{0\} \times [-1,1]$. Any tangle diagram unambiguously gives rise to a tangle, up to equivalence.
\end{Definition}

\begin{figure}
\centerline{\relabelbox 
\epsfysize 2cm 
\epsfbox{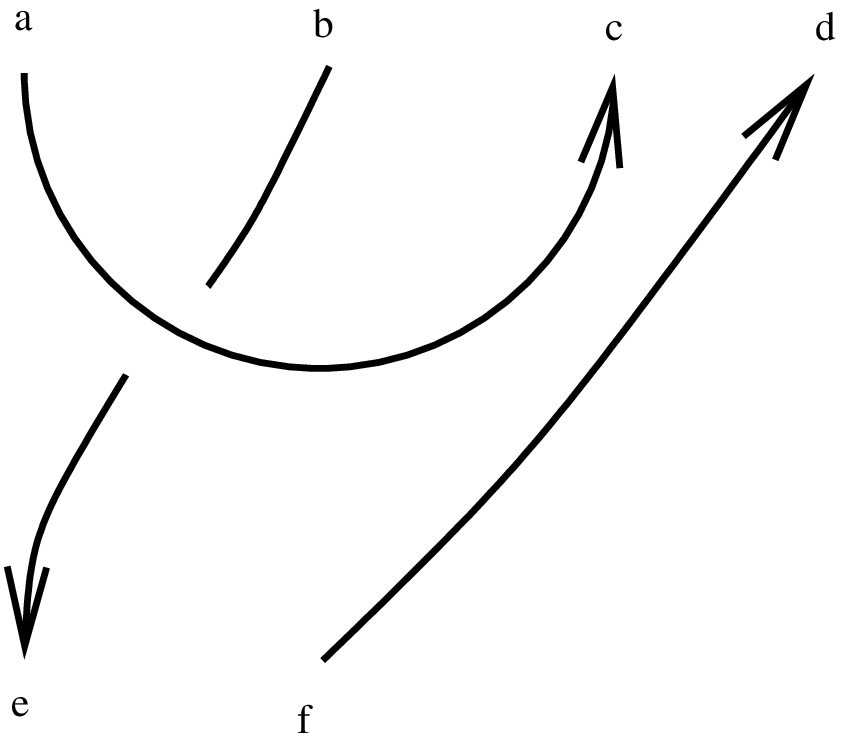}
\endrelabelbox}
\caption{A tangle with source  $++--$ with with target $+-$.\label{tangle}
}
\end{figure}

We have  monoidal categories of oriented tangles and of framed oriented tangles \cite{K,Tu}, where composition is the obvious vertical juxtaposition of tangles and the tensor product $T \tn T'$ is obtained {by placing $T'$ on the right hand side of $T$.} The objects of the categories of oriented tangles and of framed oriented tangles are words in the symbols  $\{+,-\}$; see figure \ref{tangle} for conventions. 

\begin{figure}
\centerline{\relabelbox 
\epsfxsize 10cm 
\epsfbox{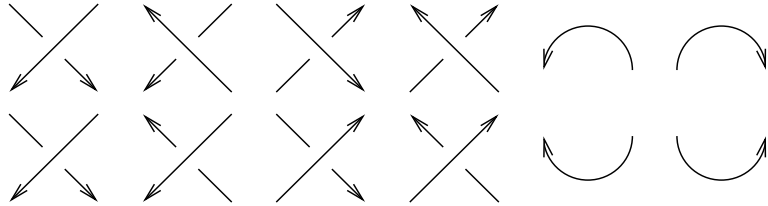}
\endrelabelbox}
\caption{Elementary generators of tangle diagrams.\label{dgenerators}}
\end{figure}

\begin{figure}
\centerline{\relabelbox 
\epsfysize 2cm 
\epsfbox{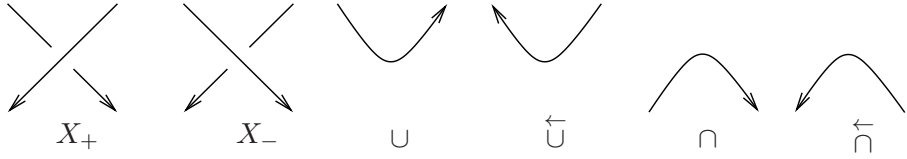}
\relabel{A}{${X_+}$}
\relabel{B}{${X_-}$}
\relabel{C}{${\cup}$}
\relabel{D}{${\stackrel{\leftarrow}{\cup}}$}
\relabel{E}{${ \cap}$}
\relabel{F}{${\stackrel{\leftarrow}{\cap}}$}
\endrelabelbox}
\caption{Generators for the categories of oriented tangles and of  framed oriented tangles.\label{tangle-generators}
}
\end{figure}

{An oriented tangle diagram is a union of the {tangle diagrams} of figure \ref{dgenerators}, with some vertical lines connecting them. }
This is a redundant set if we consider oriented tangles up to isotopy.
\begin{Definition}
 A sliced oriented tangle diagram is an oriented tangle diagram, subdivided into thin horizontal strips, inside which we have only vertical lines and possibly one of the morphisms in figure \ref{tangle-generators}.
\end{Definition}

\begin{figure}
\centerline{\relabelbox 
\epsfysize 10cm 
\epsfbox{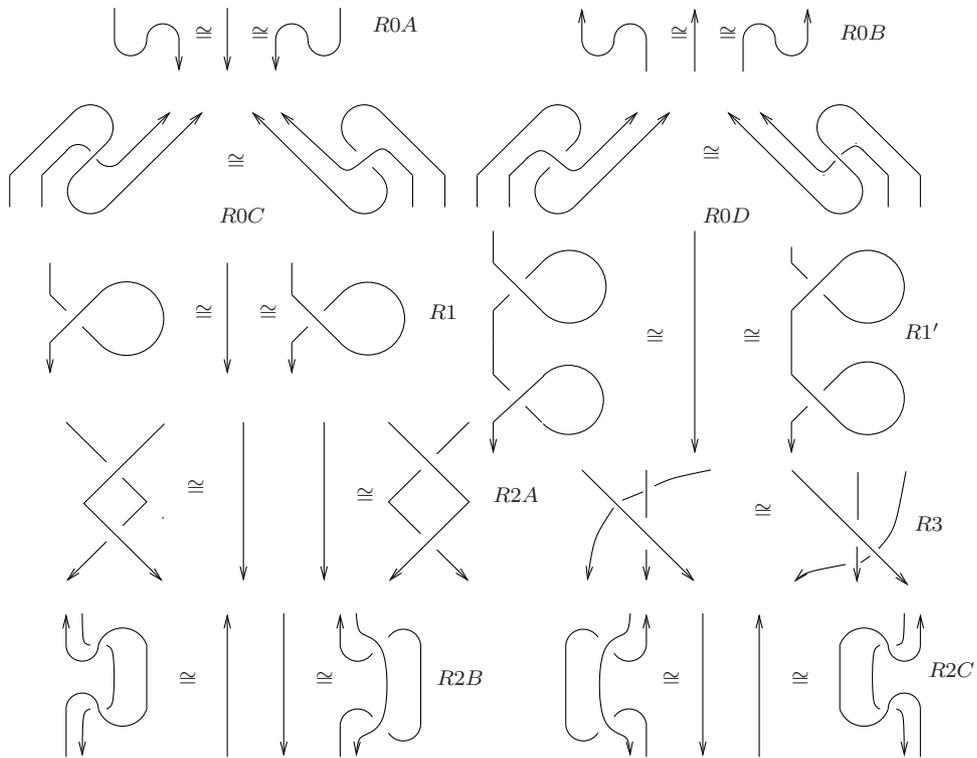}
\relabel{e1}{$\s{\cong}$}
\relabel{e2}{$\s{\cong}$}
\relabel{e3}{$\s{\cong}$}
\relabel{e4}{$\s{\cong}$}
\relabel{e5}{$\s{\cong}$}
\relabel{e6}{$\s{\cong}$}
\relabel{e7}{$\s{\cong}$}
\relabel{e8}{$\s{\cong}$}
\relabel{e9}{$\s{\cong}$}
\relabel{e10}{$\s{\cong}$}
\relabel{f}{$\s{\cong}$}
\relabel{e11}{$\s{\cong}$}
\relabel{e12}{$\s{\cong}$}
\relabel{e13}{$\s{\cong}$}
\relabel{e14}{$\s{\cong}$}
\relabel{e15}{$\s{\cong}$}
\relabel{e16}{$\s{\cong}$}
\relabel{R1}{$\s{R0A}$}
\relabel{R2}{$\s{R0B}$}
\relabel{R3}{$\s{R0C}$}
\relabel{R4}{$\s{R0D}$}
\relabel{R5}{$\s{R1}$}
\relabel{R6}{$\s{R1'}$}
\relabel{R7}{$\s{R2A}$}
\relabel{R8}{$\s{R3}$}
\relabel{R9}{$\s{R2B}$}
\relabel{R10}{$\s{R2C}$}
\endrelabelbox}
\caption{Relations for the categories of oriented tangles and of  framed oriented tangles.\label{tangle-relations}
}
\end{figure}

A theorem appearing in \cite{K}  (Theorem XII.2.2) and also in \cite{FY,RT,Tu,O} states that the category of oriented tangles {may be presented in terms of generators and relations as follows:}
\begin{Theorem}
The monoidal category of oriented tangles is equivalent to the monoidal category presented  by the six oriented tangle diagram generators :
\begin{equation}
X_+, \, X_{-},\, \cup, \, \stackrel{\leftarrow}{\cup}, \, \cap, \, \stackrel{\leftarrow}{\cap}, \,
\label{tanglegens} 
\end{equation}
shown in Figure \ref{tangle-generators}, subject to the 15 tangle diagram relations $R0A-D$, $R1$, $R2A-C$, $R3$ of Figure \ref{tangle-relations}. The category of framed oriented tangles has the same set \eqref{tanglegens} as generators, subject to the 15 relations $R0A-D$, $R1'$, $R2A-C$, $R3$ of Figure \ref{tangle-relations}. 
\end{Theorem}

\begin{figure}
\centerline{\relabelbox 
\epsfxsize 12cm 
\epsfbox{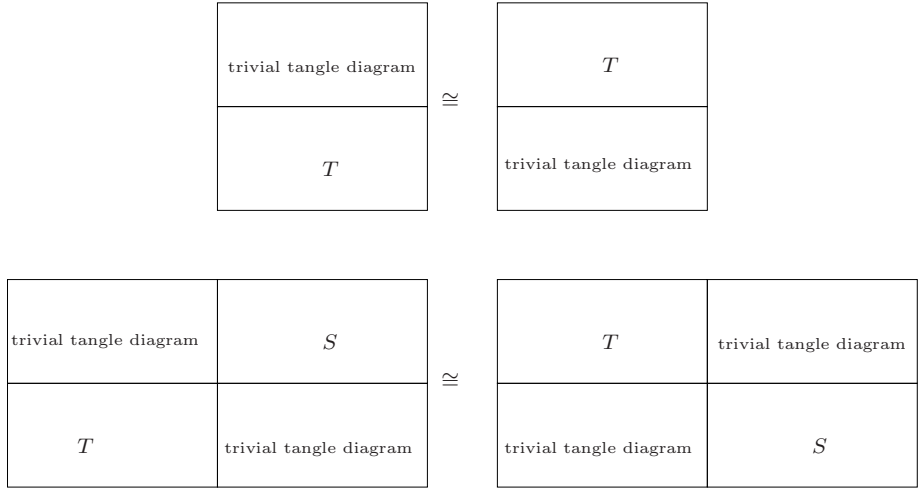}
\relabel{i1}{$\scriptscriptstyle{\textrm{trivial tangle diagram}}$}
\relabel{i2}{$\scriptscriptstyle{\textrm{trivial tangle diagram}}$}
\relabel{i3}{$\scriptscriptstyle{\textrm{trivial tangle diagram}}$}
\relabel{i4}{$\scriptscriptstyle{\textrm{trivial tangle diagram}}$}
\relabel{i5}{$\scriptscriptstyle{\textrm{trivial tangle diagram}}$}
\relabel{i6}{$\scriptscriptstyle{\textrm{trivial tangle diagram}}$}
\relabel{A}{$\s{T}$}
\relabel{B}{$\s{T}$}
\relabel{T1}{$\s{T}$}
\relabel{S1}{$\s{S}$}
\relabel{T2}{$\s{T}$}
\relabel{S2}{$\s{S}$}
\relabel{s}{$\s{\cong}$}
\relabel{t}{$\s{\cong}$}
\endrelabelbox}
\caption{The identity move and the interchange move.\label{tensor}}
\end{figure}

\begin{Remark}
We have replaced the R3 relation of Kassel's theorem with its inverse,  since this is slightly more convenient algebraically. The two forms of R3 are equivalent because of the R2A relations. The six other types of oriented crossing
in figure \ref{dgenerators}, with one or both arcs pointing upwards, can be expressed in terms of the generators (\ref{tanglegens}) and are therefore not independent generators  - see \cite{K}, Lemma XII.3.1.
\end{Remark}
The previous theorem gives generators and relations at the level of tensor categories. If we want to express not-necessarily-functorial invariants of tangles it is more useful to work with sliced tangle diagrams. The following appears for example in \cite{O}
\begin{Theorem}\label{sliced tangles}
 Two sliced oriented tangle diagrams represent the same oriented tangle (framed oriented tangle) if, and only if, they are related by 
\begin{enumerate}
 \item Level preserving isotopy of tangle diagrams.
\item The moves $R0A-R0D$, $R1$, $R2A-R2C$, $R3$ of Figure \ref{tangle-relations} (in the case of tangles), performed locally in a diagram, or the moves $R0A-R0D$, $R1'$, $R2A-R2C$, $R3$ of Figure \ref{tangle-relations}, in the case of framed tangles.
\item The ``identity'' and ``interchange'' moves of figure \ref{tensor}. (Here $T$ and $S$ can be any tangle diagrams and a trivial tangle diagram is a diagram made only of vertical lines.)
\end{enumerate}
\end{Theorem}

\begin{Definition}
[(Enhanced tangle)]\label{ent} Let $X$ be a set (normally $X$ will be either a group or a {quandle / rack}). An $X$-enhanced (framed) oriented tangle is a (framed) oriented tangle  $T$ together with an assignment of an element of $X$ to each point of the boundary of $T$. We will  consider $X$-enhanced (framed)  tangles up to isotopy of $\R \times \R \times [-1,1]$, fixing the end-points. 
\end{Definition}

\begin{figure}
\centerline{\relabelbox 
\epsfysize 2cm 
\epsfbox{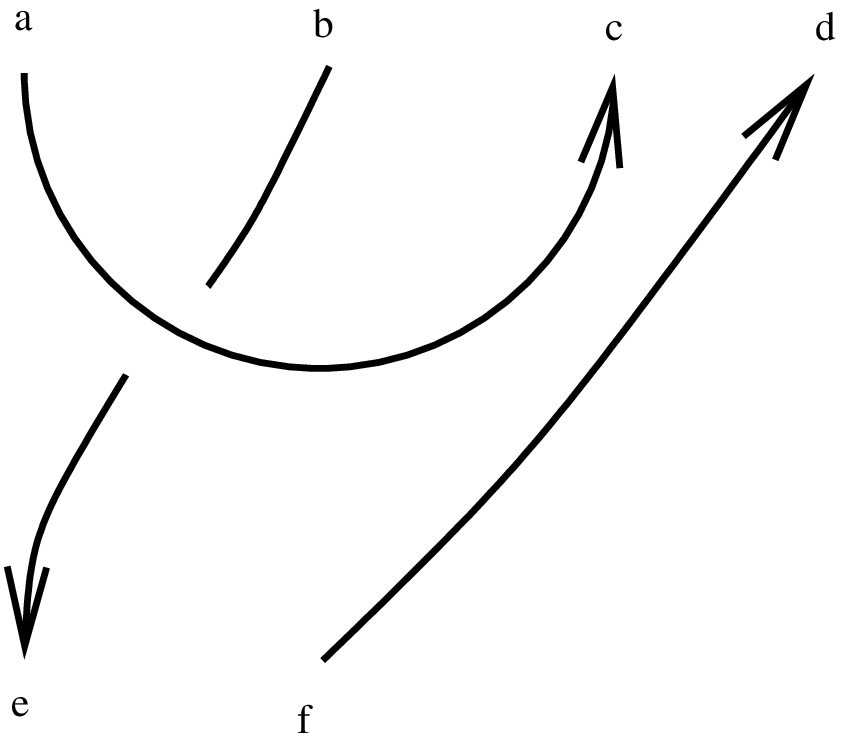}
\relabel{a}{$\s{a}$}
\relabel{b}{$\s{b}$}
\relabel{c}{$\s{c}$}
\relabel{d}{$\s{d}$}
\relabel{e}{$\s{e}$}
\relabel{f}{$\s{f}$}
\endrelabelbox}
\caption{An $X$-enhanced tangle with source $a.b.c^*.d^*$ and target $e.f^*$.\label{etangle}}
\end{figure}

Given a set $X$, there exist monoidal categories having as morphisms the set of $X$-enhanced oriented tangles and of $X$-enhanced framed oriented tangles, up to isotopy.  These categories have as objects the set of all formal words $\w$ in the symbols $a$ and $a^*$ where $a \in X$.  See figure \ref{etangle} for conventions.
\begin{Definition}\label{ev1}
Let $G$ be a group. There exists an evaluation map $\w \mapsto e(\w)$, which associates to a word $\w$ in $G \sqcup G^*$ an element of $G$, obtained by multiplying all elements of $\w$ in the same order, by putting $g^*\doteq g^{-1}$ (and $e(\emptyset)=1_G$ for the empty word {$\emptyset$}).
\end{Definition}
 \subsection{Colourings of tangle diagrams} Let ${\Gc= ( \d\colon E \to  G,\tr)}$ be a crossed module of groups.
We wish to define the notion of a $\Gc$-colouring of an oriented tangle diagram, by assigning elements of $G$ to the arcs and elements of $E$ to the crossings in a suitable way. 

For a link diagram $D$ realized as a tangle diagram, a 
$\Gc$-colouring may be regarded as a morphism of crossed modules from the fundamental crossed module $\Pi_2(X_D,Y_D)$ of {Example \ref{fundexp} and  Subsection \ref{elc}}, to $\Gc$. This idea extends in a natural way to general tangle diagrams. 

\begin{Definition}
Given a finite crossed module ${\Gc= ( \d\colon E \to  G,\tr)}$ and an oriented tangle diagram $D$, a $\Gc$-colouring of $D$ is an assignment of an element of $G$ to each arc of the diagram, and of an element of $E$ to each crossing of the diagram, such that, at each crossing of type $X_+$ or $X_{-}$ with colourings as in \eqref{GC-col}, the following relations hold:
\begin{eqnarray}
X_+: \quad \d(e) & = & XYX^{-1}Z^{-1} \label{fundX+}\\
X_-: \quad \d(e) & = & YXZ^{-1}X^{-1} \label{fundX-}
\end{eqnarray}
\end{Definition}
\begin{equation}\label{GC-col}
 \xymatrix{ &Z\ar[dr] & &X \ar[ddll]  \\
            & &{\quad\quad e}\ar[dr] &\\
            &X & &Y}  \xymatrix{ &X\ar[ddrr] & &Z \ar[dl]  \\
            & &{\quad\quad e}\ar[dl] &\\
            &Y & &X}
\end{equation}

Thus we are assigning to each type of coloured crossing a morphism of $\mathcal{C}(\Gc)$ and of $\Cl(\Gc)$, and in a similar way we may associate morphisms of $\mathcal{C}(\Gc)$ to all elementary $\Gc$-coloured tangles, as summarised in figure
\ref{tanglemor}. With the duality where the dual of the morphism $X \ra{e} \partial(e) X$ is $X^{-1} \partial(e)^{-1} \ra{X^{-1} \tr e}  X^{-1}$, and the morphisms associated to the cups and caps are the ones in figure \ref{tanglemor}, we can easily see \cite{FM} that the  monoidal category $\C(\Gc)$ is a compact category \cite{RT}, and in fact a pivotal category \cite{BW}, which however is not spherical in general. Therefore planar graphs coloured in $\C(\Gc)$ can be evaluated to give morphisms in $\C(\Gc)$, and this evaluation is invariant under planar isotopy.

Thus we can assign to the complete $\Gc$-coloured oriented tangle diagram a morphism of $\mathcal{C}(\Gc)$, by using the monoidal product horizontally and composition vertically. This leads to the following definition:
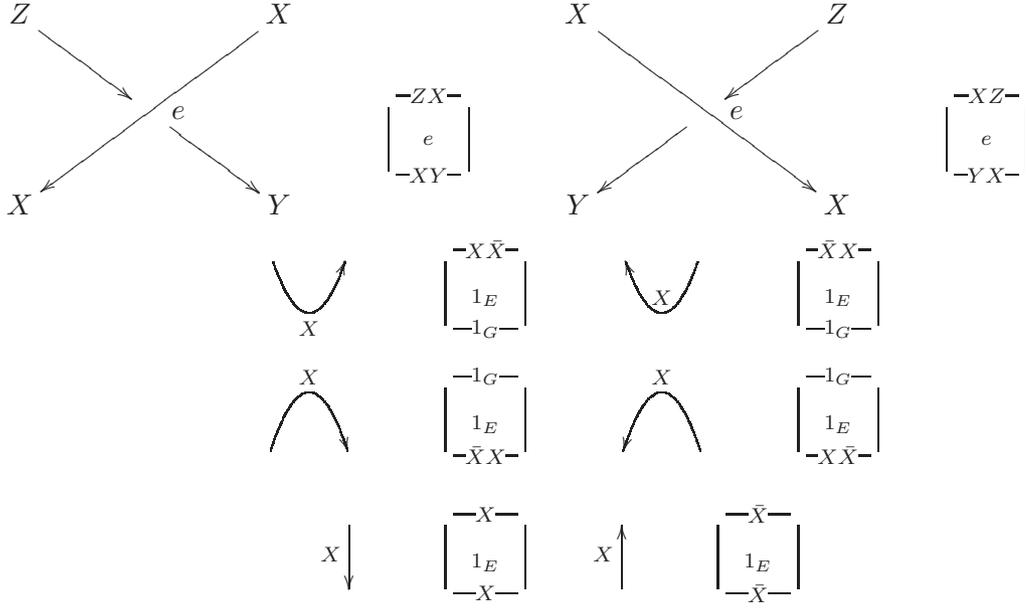
\begin{figure}
 $$\hskip-1.5cm {\xymatrix{ &Z\ar[dr] & &X \ar[ddll]  \\
            & &{\quad\quad e}\ar[dr] &\\
            &X & &Y}   \xymatrix{\\ & \ar@{-}[r]|{ZX}\ar@{-}[d] \ar@{-}@/_1pc/ @{{}{ }{}} [r]_{e}  
 & \ar@{-}[d] \\
                                                            & \ar@{-}[r]|{XY} &}  \xymatrix{ &X\ar[ddrr] & &Z \ar[dl]  \\
            & &{\quad\quad e}\ar[dl] &\\
            &Y & &X}  \xymatrix{ \\ & \ar@{-}[r]|{XZ}\ar@{-}[d] \ar@{-}@/_1pc/ @{{}{ }{}} [r]_{e}  
 & \ar@{-}[d] \\
                                                      & \ar@{-}[r]|{YX} &}}$$
$$
\xymatrix{ & \ar@/_2pc/[r]_X & } \xymatrix{  & \ar@{-}[r]|{X\bar{X}}\ar@{-}[d] \ar@{-}@/_1pc/ @{{}{ }{}} [r]_{1_E}  
 & \ar@{-}[d]\\
                                                      & \ar@{-}[r]|{1_G} &} \xymatrix{  & & \ar@/^2pc/[l]_X  } \xymatrix{ 		 & \ar@{-}[r]|{\bar{X}X}\ar@{-}[d] \ar@{-}@/_1pc/ @{{}{ }{}} [r]_{1_E}  
 & \ar@{-}[d]\\
                                                      & \ar@{-}[r]|{1_G} &} $$
$$
\xymatrix{ 		\\& \ar@/^2pc/[r]^X & } \xymatrix{  & \ar@{-}[r]|{1_G}\ar@{-}[d] \ar@{-}@/_1pc/ @{{}{ }{}} [r]_{1_E}  
 & \ar@{-}[d]\\
                                                      & \ar@{-}[r]|{\bar X X} &} \xymatrix{ \\ & & \ar@/_2pc/[l]_X  } \xymatrix{ 		 & \ar@{-}[r]|{1_G}\ar@{-}[d] \ar@{-}@/_1pc/ @{{}{ }{}} [r]_{1_E}  
 & \ar@{-}[d]\\
                                                      & \ar@{-}[r]|{X \bar X} &} $$

$$
\xymatrix{& \ar[d]_X \\& } \xymatrix{  & \ar@{-}[r]|{X}\ar@{-}[d] \ar@{-}@/_1pc/ @{{}{ }{}} [r]_{1_E}  
 & \ar@{-}[d]\\
                                                      & \ar@{-}[r]|{ X} &} \xymatrix{  & \\ & \ar[u]^X  } \xymatrix{ 		 & \ar@{-}[r]|{\bar X}\ar@{-}[d] \ar@{-}@/_1pc/ @{{}{ }{}} [r]_{1_E}  
 & \ar@{-}[d]\\
                                                      & \ar@{-}[r]|{\bar{X}} &} $$
\caption{Turning $\Gc$-coloured tangles into morphisms of $\C(\Gc)$  The symbol $X \bar{X}$ stands for $X^{-1}$.\label{tanglemor} }
\end{figure}

\begin{Definition}\label{ev2}
Given a $\Gc$-colouring $F$ of a tangle diagram $D$ (we say $F\in C_\Gc(D)$, the set of $\Gc$ colourings of $D$), the evaluation of $F$, denoted $e(F)$, is the morphism in $\mathcal{C}(\Gc)$ 
 obtained by multiplying horizontally and composing vertically the morphisms of $\mathcal{C}(\Gc)$ associated to the elementary tangles which make up $D$. 
\end{Definition}

\begin{Remark}
For a link diagram the evaluation of $F$ takes values in $A$, the automorphism subgroup of $1_G$, i.e. $A={\rm ker}\, \d\subset E$. 
\end{Remark}

For a tangle diagram without open ends at the top and bottom, i.e. a link diagram, one can conjecture that the number of 
$\Gc$-colourings of the diagram can be normalized to a link invariant, by analogy with the familiar link invariant which is the number of Wirtinger colourings of the diagram using a finite group $G$ (a Wirtinger colouring in the present context would be a $\Gc$-colouring where the group $E$ is trivial).   Indeed it was proven in \cite{FM} that the number of colourings of a link diagram evaluating to the identity of $E$ can be normalised to give an invariant of knots; see subsection \ref{elc}. However this invariant depends only on the homotopy type of the complement of the knot \cite{FM2}, thus it is a function of the knot group only.

Therefore we are led to consider the possibility of imposing more refined constraints on the $\Gc$-colourings of a tangle diagram in such a way that the number of constrained $\Gc$-colourings does respect the Reidemeister moves. Intuitively we are looking at the simple homotopy type, rather than the homotopy type of {a link complement.} Our idea is to restrict ourselves to $\Gc$-colourings of diagrams where, at each crossing, the colouring of the crossing with an element of $E$ is determined by the $G$-colouring of two arcs, namely the overcrossing arc and the lower undercrossing arc. To this end we introduce two functions:
$$
\psi: G\times G \rightarrow E, \qquad \phi: G\times G \rightarrow E,
$$
which determine the $E$-colouring of the two types of crossing, as in \eqref{psiphidef}:
\begin{equation}\label{psiphidef}
 \xymatrix{ &Z\ar[dr] \ar@{-}@/_2pc/ @{{}{ }{}}[rr]_{\quad \quad \quad \quad \quad \psi(X,Y)} & &X \ar[ddll]  \\
            & &\ar[dr] &\\
            &X & &Y}  \xymatrix{ &X\ar[ddrr] \ar@/_2pc/ @{{}{ }{}}[rr]_{\quad \quad \quad \quad \quad \phi(X,Y)}  & &Z \ar[dl]  \\
            & & \ar[dl] &\\
            &Y & &X}
\end{equation}
{Since this is  a $\Gc$-colouring, these functions  determine the $G$-assignment for the remaining arc:}
\begin{eqnarray}
X_+: \quad Z & = & \d \psi(X,Y)^{-1}XYX^{-1} \label{Zforpsi} \\
X_-: \quad Z & = & X^{-1}\d  \phi(X,Y)^{-1}YX \label{Zforphi}
\end{eqnarray}

We now come to our main definition.
\begin{Definition}[(unframed Reidemeister pair)]\label{rp}
 The pair of functions $\Phi=(\psi,\phi)$ is said to be an unframed Reidemeister pair if $\psi \colon G \times G \to E$ and $\phi\colon G \times G \to E$ satisfy the following three relations for each $ X,\, Y, \, T \in G$:
\begin{eqnarray}
 \psi(X,X) & \stackrel{R1}{=} & 1_E \label{R1}\\
 \phi(X,Y) \psi(X,Z) & \stackrel{R2}{=} & 1_E \label{R2}\\				
\f(Y,X)\,.\, Y\tr \f(T,Z)\,.\, \f(T,Y) & \stackrel{R3}{=} & X\tr \f(T,Y)\,.\, \f (T,X) \,.\, T\tr \f(V,W) \label{R3}
\end{eqnarray}
where in R2 
$$ Z = X^{-1}\d  \phi(X,Y)^{-1}YX $$
and in R3
\begin{eqnarray*}
 Z & = & Y^{-1}\d  \phi(Y,X)^{-1}XY \\
V & = & T^{-1}\d  \phi(T,Y)^{-1}YT \\ 
W & = & T^{-1}\d  \phi(T,X)^{-1}XT 
\end{eqnarray*}
\end{Definition}
\begin{Remark}
The equations above relate to the Reidemeister 1-3 moves, as we will see shortly in the proof of Theorem \ref{maintheorem}.
 If equation \eqref{R2} holds we can substitute \eqref{R3} by the equivalent:
 \begin{equation}\label{R3p} 
 \psi(X,Y)\,.\, A \tr \psi(X,Z) \,.\, \psi(A,B)= X \tr \psi(Y,Z) \,.\, \psi(X,C)\,.\, D \tr \psi(X,Y)
\end{equation}
where
\begin{equation}
 \begin{split}
  A&=\d(\psi(X,Y))^{-1} XYX^{-1}\\
  B&=\d(\psi(X,Z))^{-1} XZX^{-1}\\
  C&=\d(\psi(Y,Z))^{-1} YZY^{-1}\\
  D&=\d(\psi(X,C))^{-1} XCX^{-1}
 \end{split}
\end{equation}

\end{Remark}

\begin{figure}
\centerline{\relabelbox 
\epsfysize 3cm 
\epsfbox{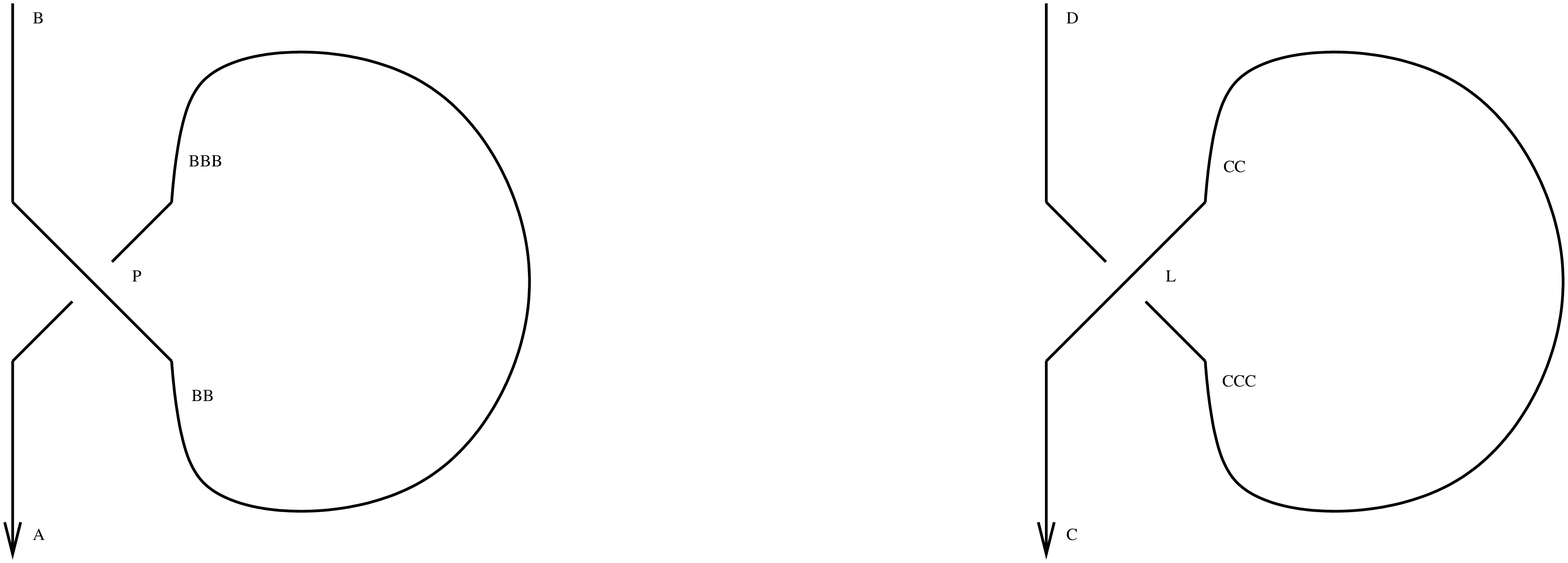}
\relabel{A}{$\scriptstyle{Z}$}
\relabel{B}{$\scriptstyle{f(Z)}$}
\relabel{BB}{$\scriptstyle{f(Z)}$}
\relabel{BBB}{$\scriptstyle{f(Z)}$}
\relabel{P}{$\scriptstyle{\phi(f(Z),Z)}$}
\relabel{C}{$\scriptstyle{A}$}
\relabel{CC}{$\scriptstyle{A}$}
\relabel{CCC}{$\scriptstyle{A}$}
\relabel{D}{$\scriptstyle{g(A)}$}
\relabel{L}{$\scriptstyle{\psi(A,A)}$}
\endrelabelbox}
\caption{Definition of $f,g \colon G \to G$.}
\label{fandg}
\end{figure}

\begin{Definition}[(Framed Reidemeister pair)]\label{frp}
 The pair of functions $\Phi=(\psi,\phi)$ is said to be a framed Reidemeister pair if relations $R2$ and $R3$ of Definition \ref{rp} hold and moreover:
\begin{enumerate}
 \item Given $Z$ in $G$, the equation (c.f. left of figure \ref{fandg}): $$\partial(\phi(A,Z))A=Z$$ has a unique solution  $f(Z) \in G$.
\item Defining $g(A)=\d(\psi(A,A))^{-1}A$ (c.f. right of figure \ref{fandg}) it holds that $f\circ g=g\circ f=\id$. In particular both $f$ and $g$ are bijective.
\end{enumerate}
\end{Definition}

\begin{Definition}
Given a finite crossed module ${\Gc= ( \d\colon E \to  G,\tr)}$, provided with a (framed or unframed) Reidemeister pair $\Phi=(\psi,\phi)$, and an oriented $G$-enhanced tangle diagram $D$, a Reidemeister $\Gc$-colouring of $D$ is a $\Gc$-colouring of $D$ (which extends the colourings at the end-points of $D$, in the sense that an arc coloured by $g$ corresponds to endpoints coloured by $g$ or $g^\ast$, depending on the orientation - see Figure \ref{etangle}) determined by the functions $\psi: G\times G \rightarrow E,\,  \phi: G\times G \rightarrow E$, which fix the colourings at each crossing as in \eqref{psiphidef}, \eqref{Zforpsi} and \eqref{Zforphi}.
\end{Definition}

We are now in a position to define a state-sum coming from the Reidemeister $\Gc$-colourings of a link diagram $D$. {Recall Definitions \ref{ev1} and \ref{ev2}.}

\begin{Definition}
Consider a  finite crossed module ${\Gc= ( \d\colon E \to  G,\tr)}$, provided with a (framed or unframed) Reidemeister pair $\Phi=(\p,\f)$. Consider an oriented enhanced tangle diagram $D$, connecting the {words $\w$ and $\w'$ in $G \sqcup G^*$}.
We denote the corresponding set of Reidemeister  $\Gc$-colourings of $D$ by $C_\Phi(D,\w,\w')$. Then we define the state sum:
\begin{equation}
I_\Phi(D) =\langle \w|I_\Phi(D)| \w' \rangle\doteq\sum_{F\in C_\Phi(D,\w,\w')} e(F)
\label{Iphi}
\end{equation}
taking values in $${\mathbb N}\big[{\rm Hom}_{\C(\Gc)}\big(e(\w),e(\w')\big)\big]\subset {\rm Hom}_{\C_{\Z}(\Gc)}\big(e(\w),e(\w')\big).$$
(Here given a category $\C$ the set of morphisms $x \to y$ is denoted by ${\rm Hom}_\C(x,y)$.)
\end{Definition}
\begin{Remark} If $D$ is a link diagram then $I_\Phi(D)$ takes values in $\Z[A]$,  the group algebra of $A = \ker \partial.$
\end{Remark}
\begin{Theorem}
The state sum $I_\Phi$ defines an invariant of $G$-enhanced tangles if $\Phi$ is an unframed Reidemeister pair and an invariant of framed $G$-enhanced tangles if $\Phi$ is a framed Reidemeister pair.
\label{maintheorem}
\end{Theorem}

\begin{Proof}
We need to show that $I_\Phi$ respects the relations of Theorem \ref{sliced tangles}. Invariance under level preserving isotopy is obvious. Let us now address, for the unframed case, the moves $R0A-D$, $R1$, $R2A-C$, $R3$ of Figure \ref{tangle-relations}. For each relation we fix matching colours on the maximum number of arcs connecting to the exterior, and then show that the corresponding morphisms of $\mathcal{C}(\Gc)$ are equal (and in some cases, that the remaining arcs connecting to the exterior are also coloured compatibly). Thus for a pair of diagrams related by one of the relations, each term in the expression for $I_\Phi$ for one diagram has a corresponding term in the expression for $I_\Phi$ for the other diagram, and the evaluations are equal term by term. 

For each relation we fix arc assignments in $G$ (see Figure \ref{tangle-relations-col}), give the corresponding equation in $E$ and show that it follows from the conditions R1-3 of Definition \ref{rp}.

\begin{figure}
\centerline{\relabelbox 
\epsfysize 10cm 
\epsfbox{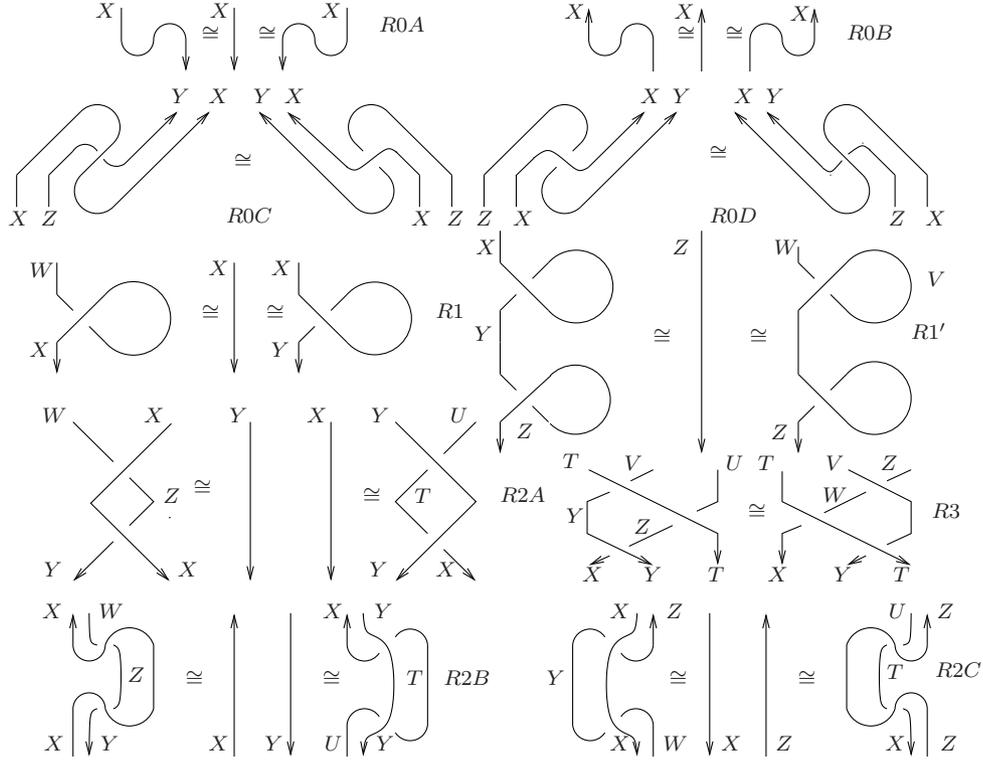}
\relabel{e1}{$\s{\cong}$}
\relabel{e2}{$\s{\cong}$}
\relabel{e3}{$\s{\cong}$}
\relabel{e4}{$\s{\cong}$}
\relabel{e5}{$\s{\cong}$}
\relabel{e6}{$\s{\cong}$}
\relabel{e7}{$\s{\cong}$}
\relabel{e8}{$\s{\cong}$}
\relabel{e9}{$\s{\cong}$}
\relabel{e10}{$\s{\cong}$}
\relabel{f}{$\s{\cong}$}
\relabel{e11}{$\s{\cong}$}
\relabel{e12}{$\s{\cong}$}
\relabel{e13}{$\s{\cong}$}
\relabel{e14}{$\s{\cong}$}
\relabel{e15}{$\s{\cong}$}
\relabel{e16}{$\s{\cong}$}
\relabel{R1}{$\s{R0A}$}
\relabel{R2}{$\s{R0B}$}
\relabel{R3}{$\s{R0C}$}
\relabel{R4}{$\s{R0D}$}
\relabel{R5}{$\s{R1}$}
\relabel{R6}{$\s{R1'}$}
\relabel{R7}{$\s{R2A}$}
\relabel{R8}{$\s{R3}$}
\relabel{R9}{$\s{R2B}$}
\relabel{R10}{$\s{R2C}$}
\relabel{x1}{$\s{X}$}
\relabel{x2}{$\s{X}$}
\relabel{x3}{$\s{X}$}
\relabel{x4}{$\s{X}$}
\relabel{x5}{$\s{X}$}
\relabel{x6}{$\s{X}$}
\relabel{x7}{$\s{Y}$}
\relabel{x8}{$\s{X}$}
\relabel{x9}{$\s{X}$}
\relabel{x10}{$\s{Z}$}
\relabel{y1}{$\s{Y}$}
\relabel{y2}{$\s{X}$}
\relabel{y3}{$\s{X}$}
\relabel{y4}{$\s{Z}$}
\relabel{y5}{$\s{Z}$}
\relabel{y6}{$\s{X}$}
\relabel{y7}{$\s{X}$}
\relabel{y8}{$\s{Y}$}
\relabel{y9}{$\s{Z}$}
\relabel{y10}{$\s{X}$}
\relabel{y11}{$\s{X}$}
\relabel{y12}{$\s{Y}$}
\relabel{z1}{$\s{W}$}
\relabel{z2}{$\s{X}$}
\relabel{z3}{$\s{X}$}
\relabel{z4}{$\s{X}$}
\relabel{z5}{$\s{Y}$}
\relabel{z6}{$\s{X}$}
\relabel{z7}{$\s{Y}$}
\relabel{z8}{$\s{Z}$}
\relabel{z9}{$\s{Z}$}
\relabel{z10}{$\s{W}$}
\relabel{z11}{$\s{V}$}
\relabel{z12}{$\s{Z}$}
\relabel{a1}{$\s{W}$}
\relabel{a2}{$\s{X}$}
\relabel{a3}{$\s{Y}$}
\relabel{a4}{$\s{Z}$}
\relabel{a5}{$\s{X}$}
\relabel{a6}{$\s{Y}$}
\relabel{a7}{$\s{X}$}
\relabel{a8}{$\s{Y}$}
\relabel{a9}{$\s{U}$}
\relabel{a10}{$\s{T}$}
\relabel{a11}{$\s{Y}$}
\relabel{a12}{$\s{X}$}
\relabel{b1}{$\s{T}$}
\relabel{b2}{$\s{V}$}
\relabel{b3}{$\s{U}$}
\relabel{b4}{$\s{Y}$}
\relabel{b5}{$\s{Z}$}
\relabel{b6}{$\s{X}$}
\relabel{b7}{$\s{Y}$}
\relabel{b8}{$\s{T}$}
\relabel{b9}{$\s{T}$}
\relabel{b10}{$\s{V}$}
\relabel{b11}{$\s{Z}$}
\relabel{b12}{$\s{W}$}
\relabel{b13}{$\s{X}$}
\relabel{b14}{$\s{Y}$}
\relabel{b15}{$\s{T}$}
\relabel{c1}{$\s{X}$}
\relabel{c2}{$\s{Y}$}
\relabel{c3}{$\s{Z}$}
\relabel{c4}{$\s{X}$}
\relabel{c5}{$\s{W}$}
\relabel{c6}{$\s{X}$}
\relabel{c7}{$\s{Y}$}
\relabel{c8}{$\s{U}$}
\relabel{c9}{$\s{Y}$}
\relabel{c10}{$\s{T}$}
\relabel{c11}{$\s{X}$}
\relabel{c12}{$\s{Y}$}
\relabel{d1}{$\s{X}$}
\relabel{d2}{$\s{W}$}
\relabel{d3}{$\s{Y}$}
\relabel{d4}{$\s{X}$}
\relabel{d5}{$\s{Z}$}
\relabel{d6}{$\s{X}$}
\relabel{d7}{$\s{Z}$}
\relabel{d9}{$\s{X}$}
\relabel{d10}{$\s{Z}$}
\relabel{d11}{$\s{T}$}
\relabel{d12}{$\s{U}$}
\relabel{d13}{$\s{Z}$}
\endrelabelbox}
\caption{Tangle relations with $G$-assignments to the arcs.\label{tangle-relations-col}
}
\end{figure}

\vskip 0.3cm
\noindent R0A and R0B. Fix $X\in G$. The corresponding equation in $E$ is $1_E=1_E$ in each case.

\vskip 0.3cm
\noindent R0C. Fix $X,Y\in G$. The corresponding equation in $E$:
$$
(X^{-1}Z^{-1})\tr \psi(X,Y) = (Y^{-1} X^{-1}) \tr \psi(X,Y)
$$
is an identity which follows from:
$$
\psi(X,Y) = (\d \psi(X,Y))\tr \psi(X,Y) = (XYX^{-1}Z^{-1}) \tr \psi(X,Y).
$$
\vskip 0.3cm
\noindent R0D. Fix $X,Y\in G$. The corresponding equation in $E$:
$$
(Z^{-1}X^{-1})\tr \phi(X,Y) = (X^{-1} Y^{-1}) \tr \phi(X,Y)
$$
is an identity which follows from:
$$
\phi(X,Y) = (\d \phi(X,Y))\tr \phi(X,Y) = (YXZ^{-1}X^{-1}) \tr \phi(X,Y).
$$

\vskip 0.3cm
\noindent R2A. Fix $X,Y\in G$. The corresponding equation in $E$ is:
\begin{equation}
\phi(X,Y) \psi(X,Z)= 1_E = \psi(Y,X)\phi(Y,T) 
\label{R2v2}
\end{equation}
where $Z=X^{-1}\d \phi(X,Y)^{-1}YX$ and $T=\d \psi(Y,X)^{-1}YXY^{-1}$.
The first equality is (\ref{R2}), and the second equality follows from the first:
$\psi(X,Z)\phi(X,Y)=1_E$ with $Z=X^{-1} \d \phi(Y,X)^{-1}XY$, i.e. $Y= \d \phi (X,Y)XZX^{-1}= 
\d \psi(X,Z)^{-1} XZX^{-1}$, then substitute variables $X\mapsto Y, \, Z\mapsto X, \, Y\mapsto T$. Applying $\d$ to (\ref{R2v2}), it follows that $W=Y$ and $U=X$.

\vskip 0.3cm
\noindent R1. Fix $X\in G$. The corresponding equation in $E$ is:
$$
\psi(X,X)=1_E = \phi(X,Y), \quad \d \phi(X,Y) =YX^{-1}
$$
The first equality is (\ref{R1}), implying $W=X$, and the second equality follows from (\ref{R1}) and (\ref{R2v2}):
$$
\phi(X,Y)=\phi(X,Y)\psi(X,X)\phi(X,X)=\phi(X,X)=1_E,
$$
which implies $Y=X$.

\vskip 0.3cm
\noindent R2B. Fix $X,Y \in G$. The corresponding equation in $E$ for the move on the left is:
$$
X^{-1}\tr \phi(X,Y)\, . \, X^{-1} \tr \psi(X,Z) = 1_E
$$
with $Z=X^{-1}\d \phi(X,Y)^{-1}YX$, which is the first equality in (\ref{R2v2}). 
Applying $\d$, it follows that $W=Y$.
The equation in $E$ for the move on the right:
$$
U^{-1} \tr \psi(Y,T) \, . \, X^{-1} \tr \phi(Y,X) = 1_E
$$
with $T=Y^{-1}\d \phi(Y,X)^{-1}XY$ and $U=\d \psi(Y,T)^{-1}YTY^{-1}$, follows from the second equality of (\ref{R2v2}):
$$
\psi(Y,T) \phi(Y,X)=1_E
$$
with $X=\d \psi(Y,T)^{-1}YTY^{-1}$, i.e.  $T=Y^{-1} \d \psi(Y,T)XY = Y^{-1}\d \phi(Y,X)^{-1}XY$, by acting with 
$X^{-1}=U^{-1}$.

\vskip 0.3cm
\noindent R2C. Fix $X,Z \in G$. The corresponding equation in $E$:
\begin{equation}
Y^{-1}\tr\phi(X,Y) \, . \, Y^{-1}\tr\psi(X,Z) = 1_E = Z^{-1}\tr \psi(Z,X) \, . \, Z^{-1}\tr \phi(Z,T)
\label{r2c}
\end{equation}
with $Y=\d\psi(X,Z)^{-1}XZX^{-1}$ and $T=\d\psi(Z,X)^{-1}ZXZ^{-1}$, is equivalent to (\ref{R2v2}) with $Y$ substituted by $Z$ in the second equality.  Applying $\d$ to (\ref{r2c}), it follows that $W=Z$ and $U=X$.

\vskip 0.3cm
\noindent R3. Fix $X,Y,T \in G$. The corresponding equation in $E$ is the Reidemeister 3 equation (\ref{R3}), which also implies the equality $U=Z$, by applying $\d$. 

\vskip 0.3cm
\noindent Invariance under the identity and interchange moves of Figure \ref{tensor} is immediate - to show the latter we use the interchange law for the operations of Figures \ref{ver} and \ref{hor}.

\vskip 0.3cm
\noindent For framed tangles we need to show invariance of $I_\Phi$ under the $R1'$ move using the properties (i) and (ii) of Definition \ref{frp}, which replace the Reidemeister 1 condition (\ref{R1}).

\vskip 0.3cm
\noindent Fix $Z\in G$. For the move on the left, at the lower crossing we have, from (i), the relation $Y=g(Z)$, and at the upper crossing we have the relation $X=f(g(Z))=Z$, by (ii).
The equation for the move:
$$
\psi(Z,Z)\, \phi(Z,Y)=1_E
$$
follows from the second equality in (\ref{R2v2}).
For the move on the right, at the lower crossing we have, from (i), the relation $V=f(Z)$, and at the upper crossing we have the relation $W=g(f(Z))=Z$, by (ii).
The equation for the move:
$$
\phi(V,Z)\, \psi(V,V)=1_E
$$
follows from the first equality in (\ref{R2v2}).

\end{Proof}

We close this section by stating a TQFT property of the invariant $I_\Phi$, which follows directly from the definitions. 

\begin{Theorem}
Let $D_1$ and $D_2$ be tangle diagrams, so that the vertical composition $\begin{array}{|l|}\hline D_1 \\ \hline D_2\\ \hline\end{array}$ is well defined. {For any enhancements $\w$ and $\w''$ of the top of $D_1$ and the bottom of $D_2$ we have:}
$$\left \langle \w \left | I_\Phi\left(\,\,\begin{array}{|l|}\hline D_1 \\ \hline D_2\\ \hline\end{array}\,\,\right) \right |\w'' \right \rangle=\sum_{\w'} \begin{array}{|l|}\hline \langle \w|I_\Phi(D_1)|\w' \rangle\\\hline \langle \w'|I_\Phi(D_2)|\w''\rangle\\ \hline \end{array} , $$
{where the sum extends over all possible enhancements $\w'$ of the intersection of $D_1$ with $D_2$.}
\end{Theorem}

\section{Examples }\label{examples}
\subsection{Examples derived from racks and quandles}\label{rackquandles}

\subsubsection{Rack and quandle invariants of links}
Recall that a rack $R$ is given by a set $R$ together with two (by the axioms not independent) operations $(x,y) \in R \times R\mapsto x\tr y \in R$  and $(x,y)\in R\times R\mapsto x \tl y \in R$, such that for each $x,y,z\in R$:
\begin{enumerate}
 \item $x\tr ( y \tl x)=y,$
\item $(x \tr y)\tl x=y,$
\item  \label{condp}$x \tr (y \tr z)=(x \tr y) \tr (x \tr z),$
\item \label{cond} $(x \tl y) \tl z=(x \tl z) \tl (y \tl z)$.
\end{enumerate}
{A quandle $Q$ is a rack satisfying moreover that:}
$$\forall x \in Q: x\tl x=x=x \tr x.$$
There is a more economical definition of a rack: it is a set $R$ with an operation $(x,y) \mapsto x \tl y$, such that condition \eqref{cond} above holds and such that for each $y \in R$ the map $ x \mapsto x \tl y$ is bijective. The inverse of it will give us the map $x \mapsto y \tr x$.

The following result and proof appeared in \cite{N}.
\begin{Lemma}[(Nelson Lemma)]\label{Nelson}
Given a rack $R$, the maps $x \mapsto x \tr x$ and $x \mapsto x \tl x$ are injective (thus bijective if the rack is finite.)
\end{Lemma}
\begin{Proof}
Let $x$ and $y$ in $R$. Then 
$$(x \tr x) \tr y=(x \tr x) \tr ( x \tr (y \tl x))=x \tr(x \tr ( y \tl x))=x \tr y.$$
If $x \tr x=y \tr y$ then
\begin{align*}
 x \tr x=(x \tr x) \tr x=( y \tr y) \tr x= y \tr x.
\end{align*}
Since $y \tr y=x \tr x$ then $y \tr y=y \tr x$. This implies $x =y$, since the map $x \mapsto y \tr x$ is bijective, its inverse being $x \mapsto x \tl y$. The proof for $x \mapsto x \tl x$ is analogous.
 \end{Proof}

\begin{figure}
\centerline{\relabelbox 
\epsfysize 2cm
\epsfbox{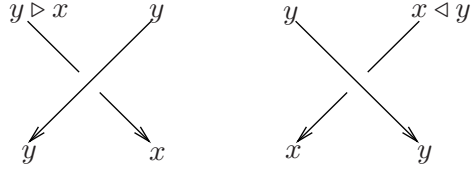}
\relabel{A}{${y \tr x}$}
\relabel{B}{${y}$}
\relabel{C}{${y}$}
\relabel{D}{${x }$}
\relabel{E}{${y}$}
\relabel{F}{${x \tl y}$}
\relabel{G}{${x}$}
\relabel{H}{${y}$}
\endrelabelbox}
\caption{A rack colouring of a link diagram in the vicinity of a vertex.
}
\label{rackcolouring}
\end{figure}

Given a knot diagram $D$, and a rack $R$, a rack colouring of $D$ is an assignment of an element of $R$ to each arc of $D$, which at each crossing of the projection has the form shown in figure \ref{rackcolouring}.

The following is well known. 
\begin{Theorem}\label{rli}
 Let $R$ be a finite rack. Then the number $I_R(D)$ of rack colourings of a link diagram $D$ is invariant under the Reidemeister moves 1', 2 and 3, and is therefore an invariant of framed links, also called $I_R$. Moreover if $R$ is a quandle then the number of rack colourings is invariant under the Reidemeister 1 move, therefore defining a link invariant.
\end{Theorem}

For a proof see \cite{N} or \cite{FR}. (Invariance under the Reidemeister moves 2 and 3 is immediate. Invariance under the Reidemeister 1' move, which is less trivial, is a consequence of the Nelson Lemma.)

\subsubsection{Reidemeister pairs derived from racks and quandles}\label{rpdrq}
Let us see that rack and quandle invariants can be written in the framework of this article. Let $R$ be a rack, which we suppose to be finite. Consider an arbitrary group structure on $R$. Call the group $G$. We do not impose any compatibility relation with the rack operations, we just assume that the underlying set of $R$ coincides with the underlying set of $G$. Consider the crossed module {$\Gc= (\id\colon G  \to G, {\rm ad})$,} where ${\rm ad}$ denotes the adjoint action of $G$ on $G$.  Put:
\begin{equation}
 \psi(B,A )=B\,\,  A  \,\, B^{-1} (B \tr A)^{-1}
\end{equation}
and
\begin{equation}
 \phi(B, A)= A\,\,B \,\,(A \tl B)^{-1} \,\,B^{-1}.
\end{equation}
\begin{Theorem}\label{racktophi}
 The pair $\Phi=(\p,\f)$ is a framed Reidemeister pair. Moreover $\Phi$ is an unframed Reidemeister pair if $R$ is a quandle.
\end{Theorem}
 \begin{Proof}
  In this case relation $R_{2}$ reads:
$$\phi(B, A) \psi(B,A\tl B)=1 $$
which follows tautologically. The relation $R_{3}$ reads in this case:
$$
\phi(B,A) \,.\, B \phi(C,A \tl B) B^{-1} \,.\, \phi(C,B) = 
A \phi(C,B)A^{-1} \,.\, \phi(C,A) \,.\, C \phi(B\tl C, A\tl C) C^{-1},
$$
and follows easily, from relation \eqref{condp} of the definition of a rack. 

Let us now prove  relations (i) and (ii) of Definition \ref{frp}.
Let $A \in G$. The equation $Z=\partial(\phi(A,Z))\,A$ means:
$$Z=Z\,A\, (Z \tl A)^{-1} , $$
or $Z \tl A=A$, that is $Z=A \tr A$. By the Nelson Lemma \ref{Nelson}, for each $Z$ the equation $Z=A \tr A$ has a unique solution $f(Z) \in R$. In this case $g(A)=\d(\psi(A,A))^{-1}A=A \tr A$. Thus trivially $f\circ g=g \circ f=\id_R$. 

Finally if $R$ is a quandle then $\psi(X,X)=X \, (X \tr X)^{-1}=X \, X^{-1}=1_E$.

 \end{Proof}

Since there is clearly, {by \eqref{psiphidef}, \eqref{Zforpsi} and \eqref{Zforphi}, a one-to-one correspondence between $\Gc$-colourings of a link diagram $D$ and rack colourings (with respect to $R$) of $D$, we have:}
\begin{Theorem}\label{rackex}
 Given a link diagram $D$:
$$I_\Phi(D)=I_R(D) 1_G,$$
where $1_G$ is the identity of $G$.
Therefore the class of invariants defined in this paper is at least as strong as {the class of rack link invariants.}
\end{Theorem}

There is a spin-off of the rack invariant in order to handle tangles. Given a rack $R$, recall that an $R$-enhanced tangle, Definition \ref{ent}, is a tangle together with a map from the boundary of $T$ into (the underlying set of) $R$. There is a category whose objects are the words $\w$ in $R \sqcup R^*$ and whose morphisms are $R$-enhanced tangles connecting them. Thus if the word $\w$ is the source of $T$ then $\w$ is a word having $i$ elements, where $i$ is the number of intersections of the tangle with $\R^2 \times \{1\}$. Moreover the $n^{\rm th}$ element of $\w$ is either the colour $a \in R$ given to the $n$-th intersection, or it is   $a^*$, the former happening if the strand is pointing downwards and the latter if the associated strand is pointing upwards. These conventions were explained in figure \ref{etangle}. Given an $R$-enhanced tangle $T$ let $\w(T)$ and $\w'(T)$ be the source and target of $T$, both words in $R \sqcup R^*$. Define $\langle \w(T) |I_R(T)| \w'(T) \rangle$ as being the number of  arc-colourings of a diagram of $T$ extending the enhancement of $T$. This defines an invariant of tangles, for any choice of colourings on the top and bottom of $T$ (in other words, for any $R$-enhancement of $T$). Clearly
\begin{Theorem}\label{tanglerack}
For any $R$-enhanced tangle $T$, putting $\w=\w(T)$ and $\w'=\w'(T)$ we have:
$$ \langle \w(T) |I_\Phi(T)|\w'(T)\rangle =   \langle \w |I_R(T)|\w'\rangle \quad \begin{CD} e(\w)\\@VV e(\w') e(\w)^{-1}V \\ e(\w')\end{CD}$$
\end{Theorem}
\noindent See definitions \ref{ev1} and \ref{ev2} for notation.
\subsubsection{Reidemeister pairs derived from rack and quandle cocycles}\label{qc}
We can extend the statement of Theorem \ref{rackex} for the case of rack cohomology invariants of knots. Let $R$ be a rack. Let $V$ be an abelian group. We say that a map $w \colon R \times R \to V$ is a rack 2-cocyle if:
$$w(x,y)+w(x \tl y,z)=w(x,z)+ w(x \tl z, y \tl z), \textrm{ for each}\, x,y,z \in R. $$
If $R$ is a quandle, such a $w$ is said to be a quandle cocycle if moroever $w(x,x)=0_V$, for each $x \in R$. For details see \cite{CJKLS,CJKLS2,N,E1,E2}.

Consider any group structure  $G$ on the set $R$, which may be completely independent of the 
rack operations. Consider the crossed module $(\partial \colon G \times V \to G,\bullet)$, where $$g \bullet (h,v)=(ghg^{-1},v), \textrm{ for each } g,h \in G \textrm{ and } v \in V,$$ which is a left action of $G$ on $G \times V$ by automorphisms, and $\partial(g,v)=g$, for each $(g,v) \in G \times V$. 

Given a rack 2-cocycle $w\colon R \times R \to V$, set:
\begin{equation}
 \psi(B,A )=\big (B\,\,  A  \,\, B^{-1} (B \tr A)^{-1},w(B \tr A,B)\big),
\end{equation}
and
\begin{equation}
 \phi(B, A)= \big (A\,\,B \,\,(A \tl B)^{-1} \,\,B^{-1}, w(A,B)^{-1}\big).
\end{equation}
\begin{Theorem}
 The pair $\Phi=(\psi,\phi)$ is a framed Reidemeister pair. Its associated framed link invariant coincides with the usual rack cohomology invariant of framed links. Morever if $R$ is a quandle and $w$ a quandle 2-cocycle then $\Phi=(\psi,\phi)$ is an unframed Reidemeister pair and its associated invariant  of links coincides with the usual quandle cocycle link invariants \cite{CJKLS}. 
\end{Theorem}
\begin{Proof}
Analogous to the proof of Theorems \ref{racktophi} and \ref{rackex}.
\end{Proof}

\noindent Therefore the class of invariants defined in this paper is at least as strong as the class of invariants of links derived from quandle cohomology classes.

\subsection{Relation with the Eisermann knot invariant}\label{ep}

\subsubsection{String knots, long knots, knot meridians and knot longitudes}\label{sklk}
Recall that an (oriented) long knot is an embedding  $f$ of  $\R$ into $\R^3$. such that, for sufficiently large (in absolute value) $t$, we have {$f(t)=(0,0,-t)$.} These are considered up to isotopy with compact support. Clearly long knots (up to isotopy) are in one-to-one correspondence with isotopy classes of tangles whose underlying 1-manifold is the interval, and whose boundary is $\{0\} \times \{0\} \times \{\pm 1\}$, being, furthermore, oriented downwards. These are usually called string knots.

There exists an obvious closing map $c$ sending a string knot $L$ to a closed knot ${\rm cl}(L)$. It is well known that this defines a one-to-one correspondence between isotopy classes of string knots and isotopy classes of oriented knots. To see this, note that a map sending a closed knot $K$ to a long knot $L_K$ can be obtained by choosing a base point $p \in K$. Then there exists an (essentially unique) orientation  (of $S^3$ and of $K$) preserving  diffeomorphism $(S^3\setminus \{p\},K \setminus \{p\}) \to (\R^3,L_K^p)$, where $L_K^p$ is a long knot with ${\rm cl}(L_K^p)=K$. Note that $L_K^p$ depends only on the orientation preserving diffeomorphism class of the triple $(S^3,K,p)$, thus since all pairs $(K,p)$, with fixed $K$, but arbitrary $p$, are isotopic we can see that 
$L_K^p$ depends only on $K$, thus we can write it as $L_K$.
\begin{figure}
\centerline{\relabelbox 
\epsfysize 2cm 
\epsfbox{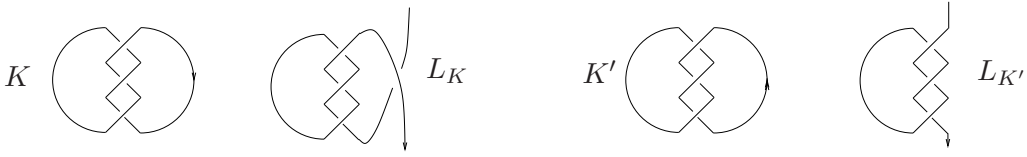}
\relabel{K}{${K}$}
\relabel{LKK}{${L_{K'}}$}
\relabel{LK}{${L_{K}}$}
\relabel{KK}{${K'}$}
\endrelabelbox}
\caption{Turning an oriented knot into a string knot in two different cases.\label{LinkTangle}
}
\end{figure}

Let $D$ be a knot diagram of the knot $K$. Consider the Wirtinger generators of the fundamental group $\pi_1(C_K)$ of the complement $C_K=S^3 \setminus n(K)$ of $K$; we thus have a meridian for any arc of the diagram\ $D$. Let $a$ be an arc of $D$ and $p$ a base point of $K$ in $a$.  Then there is a meridian $m_p=m$ of $K$ encircling {the arc $a$ at $p$}, whose direction is determined by the right hand rule. Let $D^2=\{z \in \mathbb{C}: |z|\leq 1\}$ and $S^1 = \partial D^2$. Choose an embedding $f\colon S^1\times D^2 \to S^3$ such that 
\begin{itemize}
 \item $f(S^1 \times \{0\})=K$, preserving orientations, with $f(1,0)=p$.
\item $f(\{1\} \times S^1)=m$
\item $f(S^1 \times \{1\})$ has zero linking number with $K$.
\end{itemize}
If we take  $f(1,1)$ to be the base point of $S^3$, then the homotopy class $l_p=l\in \pi_1(C_K)$  of   $f(S^1 \times \{1\})$ is called a longitude of $K$ {\cite{BZ}}. It is well known that the triple $(\pi_1(C_K), m,l)$, considered up to isomorphism, is a complete invariant of the knot $K$ \cite{W}. Note that if we choose another base point $p'$ of $K$ then $m_{p'}$ and $l_{p'}$ can be obtained from $m_p$ and $l_p$ by conjugating by a single element of $\pi_1(C_K)$.

The longitude $l_p$, being an element of the fundamental group of the complement {$C_K$ of $K$}, can certainly be expressed in terms of the Wirtinger generators. This can be done in the following way;  for details see \cite{E1,E2}. Let $a_0=a$ be the arc of the diagram $D$ of $K$ containing the base point of $K$. Then go around the knot in the direction of its orientation. This makes it possible to order the arcs of $K$, say as  $a_1,a_2,\dots, a_n$;  we would have $a_n=a_1$, except that we prefer to see $K$ as being split at the base point $p$, separating the arc $a$ in two. We can also order the crossings of $D$.

\begin{figure}
\centerline{\relabelbox 
\epsfysize 4.3cm
\epsfbox{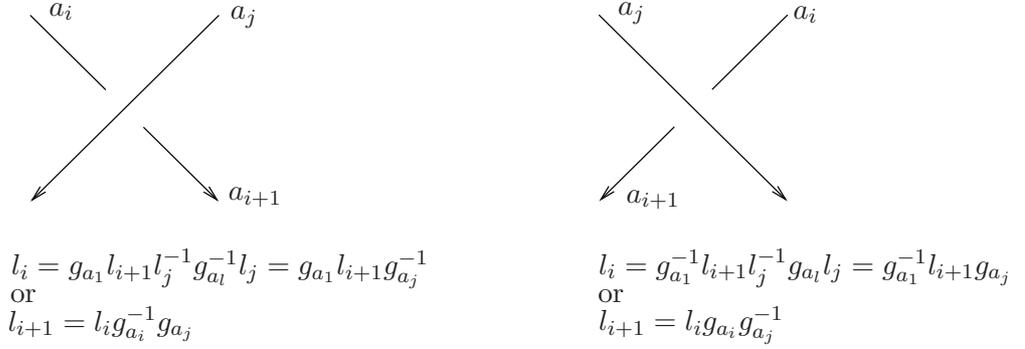}
\relabel{a}{${a_j}$}
\relabel{l}{${a_i}$}
\relabel{ll}{$a_{i+1}$}
\relabel{b}{${a_j}$}
\relabel{k}{${a_i}$}
\relabel{kk}{$a_{i+1}$}
\relabel{M}{$l_{i}=g_{a_1}l_{i+1}l_{j}^{-1}g_{a_l}^{-1}l_{j}=g_{a_1}l_{i+1}g_{a_j}^{-1}$}
\relabel{MM}{$l_{i+1}=l_i g_{a_{i}}^{-1}g_{a_j} $}
\relabel{N}{$l_{i}=g_{a_1}^{-1}l_{i+1}l_{j}^{-1}g_{a_l}l_{j}=g_{a_1}^{-1}l_{i+1}g_{a_j}$}
\relabel{NN}{$l_{i+1}=l_i g_{a_{i}}g_{a_j}^{-1} $}
\relabel{A}{or}
\relabel{P}{or}
\endrelabelbox}
\caption{Rules for partial longitudes at crossings.}
\label{partial}
\end{figure}

 The longitude $l_p$ of $K$ is expressed as a product of all elements of {$\pi_1(C_K)$, associated to the arcs we undercross} as we travel from $p$ to $p$, making sure that the linking number of $l_p$ with $K$ is zero. {Therefore any arc $a_i$ has also assigned a partial longitude $l_i$ (the product of the elements  of {$\pi_1(C_K)$, associated to the arcs we undercross,  as we travel from $p$ to $a_i$)}. We thus have $l_n=l$.} 
Given an arc  $a$ of $D$ denote the corresponding element of the fundamental group of the complement  by  $g_a$. Then clearly we have $$g_{a_i}=l_i^{-1} g_{a_1} l_i.$$  The way to pass from  $l_i$ to $l_{i+1}$ appears in figure \ref{partial} for the positive and negative crossing. 

Given $i$, let $j_i$ be the number of the arc separating $a_i$ from $a_{i+1}$. Let $\theta_i$ be the sign of the $i$-th crosssing. Then:
\begin{equation}\label{longformula}
 l=\prod_{i=1}^{n-1} g_{a_i}^{-\theta_i}g_{a_{j_i}}^{\theta_i}=\prod_{i=1}^{n-1} l_i^{-1}g_{a_1}^{-\theta_i}l_i l_{j_i}^{-1} g_{a_1}^{\theta_i} l_{j_i}=\prod_{i=1}^{n-1}[l_i^{-1},g_{a_1}^{-\theta_i}]\,\,[{g_{a_1}^{-\theta_i}},l_{j_i}^{-1}];
\end{equation}
more generally, if $k \in \{2,\dots,n\}$:
\begin{equation}\label{partiallongformula}
 l_k=\prod_{i=1}^{k-1} g_{a_i}^{-\theta_i}g_{a_{j_i}}^{\theta_i}=\prod_{i=1}^{k-1} l_i^{-1}g_{a_1}^{-\theta_i}l_i l_{j_i}^{-1} g_{a_1}^{\theta_i} l_{j_i}=\prod_{i=1}^{k-1}[l_i^{-1},g_{a_1}^{-\theta_i}]\,\,[{g_{a_1}^{-\theta_i}},l_{j_i}^{-1}].
\end{equation}
Therefore both the longitude $l$ and any partial longitude $l_k$ are elements of the commutator subgroup of the fundamental group of the complement of $K$.
Also:
$$l_{i+1}= l_i [l_i^{-1},g_{a_1}^{-\theta_i}]\,\,[{g_{a_1}^{-\theta_i}},l_{j_i}^{-1}].$$

Later on we will present another formula for a knot longitude.

\subsubsection{The Eisermann invariant of knots} Let $K$ be a knot in $S^3$. 
 Consider the fundamental group of the complement  $C_K=S^3 \setminus n(K)$ of the knot $K$. Here $n(K)$ is an open regular neighbourhood of $K$. Choose a base point $p$ of $K$. Let the associated meridian and longitude of $K$ in $\pi_1(C_K)$ be denoted by $m_p$ and $l_p$, respectively. Note that $[m_p,l_p]=1$.

Let $f\colon \pi_1(C_K) \to G$ be a group morphism. {Therefore}: $f(l_p)\in G'\doteq [G,G]$, the derived (commutator) group of $G$, generated by the commutators $[g,h]=ghg^{-1}h^{-1}$, where $g,h \in G$. Moreover {$[f(l_p),f(m_p)]=1_G$.} Then $$f(l_p)\in \Lambda\doteq [G,G] \cap C(x),$$
where {$x=f(m_p)$} and $C(x)$ is the set of elements of $G$ commuting with $x$.

Let $G$ be a finite group. Let $x$ be an element of $G$. The Eisermann invariant  \cite{E2} (also called Eisermann polynomial) is:
$$E(K)=\sum_{\left \{\substack{f\colon \pi_1(C_K) \to G \\ f(m_p)=x}\right \}} f(l_p) \in \N(\Lambda). $$
Note that if we choose a  different base point $p'$ of $K$ then $E(K)$ stays invariant since $m_{p'}=h\,m_p\,h^{-1}$ and $l_{p'}=h\,l_p\, h^{-1}$, for some common {$h \in \pi_1(C_K)$}. The Eisermann invariant can be used to detect chiral and non invertible knots  \cite{E2}.

Clearly $E(K)$ is given by a map $f_E^K\colon G' \to \N$, where 
$$E(K)=\sum_{g \in G'}  f_E^K(g) g. $$
Note that $f_E^K(g)=0$ if $g \not \in  \Lambda$.

Let us see that the Eisermann invariant can be addressed using Reidemeister pairs. This is  a consequence of 	the previous subsections and the discussion in \cite{E2,E1}, which we closely follow, having discussed and completed the most relevant issues in \ref{sklk}.

Let $G$ be a group. Choose $x\in G$ and consider from now on the pair $(G,x)$. We will use the notation $h^g=g^{-1} h g$, where $g,h \in G$.

Let \begin{align*}
Q&=\big\{x^g, g \in G'\big\} \subset G,\\
\overline{Q}&=\big\{x^g, g \in G\big\} \subset G.
\end{align*}

\begin{Lemma}
 Both sets $Q$ and $\overline{Q}$ are self conjugation invariant:
$$a,b \in Q \implies a^{-1} b a \in Q$$
and
$$a,b \in \overline{Q} \implies a^{-1} b a \in \overline{Q}$$
Therefore $Q$ and $\overline{Q}$ are both quandles, with quandle operation $h\tl g=h^g$. 
\end{Lemma}
\begin{Proof}
 Given $g,h \in G'$ we have $$(g^{-1} xg )^{-1} (h^{-1} x h)  (g^{-1} xg )=x^{x^{-1}hg^{-1} x g}$$
and $x^{-1}hg^{-1} x g=x^{-1}hg^{-1} x gh^{-1} h=[x^{-1},hg^{-1}] h \in G'$. 
The proof for $\overline{Q}$ is analogous.
\end{Proof}

It is easy to see that:
\begin{Lemma}[(Eisermann)] \label{eqd}
Let $G$ be a group. Given arbitrary $x \in G$, both $G'$ and $G$ are quandles, with quandle operation: 
$$h\tl g=x^{-1}hg^{-1} x g,$$
and therefore:
$$g \tr h'=xh'g^{-1}x^{-1} g. $$
There are also  quandle maps $p\colon G' \to Q$ and $\overline{p}\colon \overline{Q} \to G$ with $p(g)=x^g$.
\end{Lemma}
Recall the rack  invariant of tangles, defined just before Theorem \ref{tanglerack}.
\begin{Theorem}[(Eisermann)]
For any  knot $K$ and any $g$ in $G'$ it holds that:
$$\langle 1_G |I_{G'}({D_K})| g\rangle =f^K_E(g), $$
where $D_K$ is any string knot diagram associated to $K$. Of course we regard $D_K$ as an enhanced tangle diagram coloured with $1_G=1_{G'}$  at the top and with $g$ at the bottom.
\end{Theorem}
\begin{Proof}
 Follows from the discussion in \ref{sklk} and especially figure \ref{partial}.
\end{Proof}

\begin{Remark}
 The previous theorem is also valid for the quandle structure in  $G$, Lemma \ref{eqd}. Given the form of the quandle is easy to see that if $g,h \in G'$ then:
$$\langle g |I_{G'}({D_K})|h\rangle=\langle g |I_{G}({D_K})|h\rangle.$$
\end{Remark}

By using subsection \ref{rackquandles} we will show that the Eisermann invariant can be addressed in our framework, by passing to string knots. 
 Suppose we are given  a finite group $G$ and $x \in G$ {(it may be that $x \in G \setminus G'$)}.
  We can choose any group operation structure in the underlying set of $G'$. We take the most obvious one, given by the inclusion $G' \subset G$.  The associated crossed module is $G'\stackrel{\rm id}{\rightarrow} G'$, with $G'$ acting on itself by conjugation.

 {Define, given $L,M \in G'$, the pair  $\Phi^x=(\psi^x,\phi^x)$,  as:}
\begin{equation}
\begin{split}
\phi^x(L,M)&= ML(x^L)^{-1} x^M M^{-1}L^{-1}=Mx^{-1}LM^{-1}x L^{-1}=[Mx^{-1},Lx^{-1} ]\\
\psi^x(L,M)&=LML^{-1} L^{-1} x LM^{-1}x^{-1}=[L,M][ML^{-1},x]=[xML^{-1}x^{-1}Lx^{-1},Lx^{-1}]^{-1}
\end{split}
\label{eisphipsi}
\end{equation}
\begin{Theorem}
 The pair $\Phi^x=(\psi^x,\phi^x)$ is an unframed Reidemeister pair for the crossed module $G'\stackrel{\rm id}{\rightarrow} G'$, with $G'$ acting on itself by conjugation. Let $K$ be an oriented knot and $L_K$ be the associated string knot.  Given $g \in G'$ then
\begin{equation}\label{eqinv}\langle 1_{G'} | I_{\Phi^x}(L_K)| g\rangle=f^K_E(g).\end{equation}
\label{eisermann}
\end{Theorem}
\begin{Proof}
The expressions for {$\Phi^x=(\psi^x,\, \phi^x)$} guarantee that the colourings of the arcs at a crossing are those given by the Eisermann quandle operation and its inverse (Lemma \ref{eqd}). Thus  equation \eqref{eqinv} holds, and it is enough to check that $\Phi$ does indeed satisfy the conditions to be a unframed Reidemeister pair. Clearly the Reidemeister 1 condition \rref{R1} holds: $\psi(L,L)=1$.  The Reidemeister 2 equation \rref{R2} is:
$$
\phi^x(L,M)\, \psi^x(L, x^{-1}ML^{-1}xL) =1
$$  
i.e.
$$
[Mx^{-1},Lx^{-1}]\, [L,x^{-1}ML^{-1}xL ]\, [x^{-1}ML^{-1}x,x]=1.
$$
Writing this out in full, one obtains:
\begin{eqnarray*}
\lefteqn{Mx^{-1}Lx^{-1}xM^{-1}xL^{-1} \, . \, L(x^{-1}ML^{-1}xL)L^{-1}(L^{-1}x^{-1}LM^{-1}x) \, . \, }\\
& & (x^{-1}ML^{-1}x)x(x^{-1}LM^{-1}x)x^{-1} \\
 & & = M\, \underline{x^{-1}LM^{-1}xL^{-1} \, . \, L x^{-1} ML^{-1}x} \, \,\,
\underline{L^{-1}x^{-1} L M^{-1} x\, . \, x^{-1} ML^{-1} x L }\, M^{-1} = 1.
\end{eqnarray*}
In the above computation, {the underlined factors are} equal to 1.

To avoid confusion with the quandle operation $\tr$, the left action of $G'$ on $G'$ by conjugation will be denoted by $g \bullet h$, thus $g \bullet h=ghg^{-1}$ for each $g,h \in G'$. 
The Reidemeister 3 equation \rref{R3}
for this case reads:
$$
\phi^x(L,M)\, . \, L\bullet \phi^x(N,P) \, . \, \phi^x(N,L) = M\bullet \phi^x(N,L) \, . \,\phi^x(N,M) \, . \, N\bullet \phi^x(R,Q),
$$ 
where
$$
P=x^{-1}ML^{-1}xL, \quad \quad R= x^{-1}LN^{-1}xN, \quad \quad Q=x^{-1}MN^{-1}xN.
$$
The left-hand side of the above equation, written out in full, is:
\begin{eqnarray*}
\lefteqn{Mx^{-1}LM^{-1}xL^{-1} \, . \, L(Px^{-1}NP^{-1}xN^{-1})L^{-1}  \, . \, (Lx^{-1}NL^{-1} xN^{-1} )} \\
& & = M\, \underline{x^{-1}LM^{-1}x \, . \, x^{-1} ML^{-1}x } \, L x^{-1} N (L^{-1}x^{-1}LM^{-1}x)x  
N^{-1} \, . \, x^{-1} N L^{-1} x N^{-1} \\
& &  = MLN (N^{-1}x^{-1}NL^{-1}x^{-1}LM^{-1}x^2) (N^{-1}x^{-1}NL^{-1}x)N^{-1},
\end{eqnarray*}
and the right-hand side leads to the same expression:
\begin{eqnarray*}
\lefteqn{M(Lx^{-1}NL^{-1}xN^{-1})M^{-1} \, .  \,  Mx^{-1}NM^{-1}xN^{-1}  \, . \,   N(Qx^{-1}RQ^{-1}xR^{-1})N^{-1} }\\
& & = M Lx^{-1}NL^{-1}  \,  \underline{xN^{-1}   x^{-1}NM^{-1}x  (x^{-1} MN^{-1}x N) \, . \, x^{-1} } \\ 
 & & (x^{-1}L \, \underline{N^{-1} x N)( N^{-1} x^{-1} N }\, M^{-1}x)   x(N^{-1}x^{-1}NL^{-1}x)N^{-1} \\
& & =  MLN (N^{-1}x^{-1}NL^{-1}x^{-1}LM^{-1}x^2) (N^{-1}x^{-1}NL^{-1}x)N^{-1}.
\end{eqnarray*}
In both derivations we insert the definitions of $P,\, R$ and $Q$ in the first equality, and regroup the factors after eliminating the underlined expressions in the second equality. Note that both sides equal $MLNU^{-1}R^{-1}N^{-1}$, which is the product of the colourings of the 6 external arcs in Figure \ref{R3Eisermann} for Reidemeister 3, taken in anticlockwise order starting with $M$, where $U$ is the colouring assigned to the rightmost upper arc, i.e. 
$$U= x^{-1}PN^{-1}xN=x^{-1} QR^{-1}xR.$$
\begin{figure}
\centerline{\relabelbox 
\epsfysize 3cm
\epsfbox{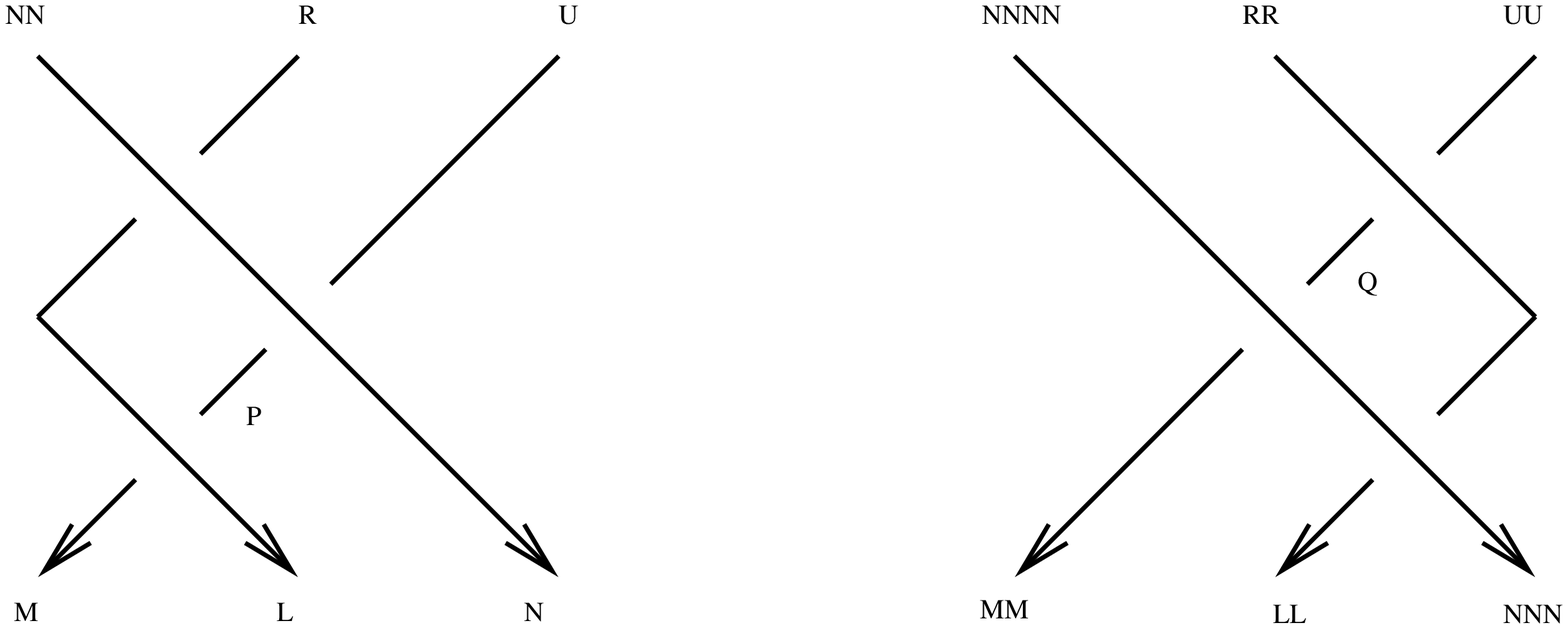}
\relabel{M}{$\scriptstyle{M}$}
\relabel{L}{$\scriptstyle{L}$}
\relabel{N}{$\scriptstyle{N}$}
\relabel{NN}{$\scriptstyle{N}$}
\relabel{R}{$\scriptstyle{R}$}
\relabel{U}{$\scriptstyle{U}$}
\relabel{P}{$\scriptstyle{P}$}
\relabel{MM}{$\scriptstyle{M}$}
\relabel{LL}{$\scriptstyle{L}$}
\relabel{NNN}{$\scriptstyle{N}$}
\relabel{NNNN}{$\scriptstyle{N}$}
\relabel{RR}{$\scriptstyle{R}$}
\relabel{UU}{$\scriptstyle{U}$}
\relabel{Q}{$\scriptstyle{Q}$}
\endrelabelbox}
\caption{Two sides of the Reidemeister-III move in the proof of Theorem \ref{eisermann}.}
\label{R3Eisermann}
\end{figure}
\end{Proof}

Consider the quandle structure in $G$, Lemma \ref{eqd}, given by the same formulae as the one of $G'$. Consider the crossed module given by the identity map $G \to G$ and the adjoint action of $G$ on $G$.  Given $x \in G$ we have an unframed Reidemeister pair  $\bar{\Phi}^x=(\bar\psi^x,\bar \phi^x)$, with the same formulae as \eqref{eisphipsi}, namely if $x \in G$ and $L,M \in G$:
 \begin{equation}\label{epb}
\begin{split}
\bar{\phi}^x(L,M) &= ML(x^L)^{-1} x^M M^{-1}L^{-1}=Mx^{-1}LM^{-1}x L^{-1}=[Mx^{-1},Lx^{-1}]\\
\bar{\psi}^x(L,M)&=LML^{-1} L^{-1} x LM^{-1}x^{-1}=[L,M][ML^{-1},x]=[xML^{-1}x^{-1}Lx^{-1},Lx^{-1}]^{-1}
\end{split}
\end{equation}
{It is easy to see that:}
\begin{Proposition}
 Let $a,b \in G'$. Then for each $x \in G$:
$$\langle a | I_{\Phi^x}(L_K)| b\rangle =\langle a  | I_{\bar{\Phi}^x} (L_K)| b\rangle$$
for each string knot $L_K$.
\end{Proposition}

\begin{Remark}[(Formula for a knot longitude)]\label{longformulanew}
Our approach for defining Eisermann invariants, and the proof of Theorem \ref{eisermann}, provides a different formula to \eqref{longformula} for the longitude {$l_p=l$} of a knot $K$, if $K$ is presented as the closure of a string knot $L$. It is assumed that the base point ${p \in K}$ of the closed knot $K$ lives in the top end of $L$.  Let $G=\pi_1(C_K)$. Consider the crossed module $(\id\colon G \to G,\ad)$. Let $x$ be the element of $G$ given by the top strand  $a$ of $L$. Let $b$ be the bottom strand {of $L$}.  Consider the Reidemeister pair $\bar{\Phi}^x$ in \eqref{epb}, thus
\begin{equation}
\bar{\phi}^x(L,M)=[Mx^{-1},Lx^{-1} ] \textrm{ and } \bar{\psi}^x(L,M)=[xML^{-1}x^{-1}Lx^{-1},Lx^{-1}]^{-1}.
\end{equation}
 Consider a diagram $D$ of $L$.
 Colour each arc $c$ of the diagram $D$  with the corresponding partial longitude $l_c$, as defined in \ref{sklk}. Therefore $l_a=1$ and $l_b=l$. Then by the proof of Theorem \ref{eisermann} one has a Reidemeister colouring $F$. If we evaluate $F$ {(definitions \ref{ev1} and  \ref{ev2})}, we have a morphism $l_a \ra{e(F)} l_b$,  hence {$e(F)=l=l_p$.} The form of $e(F)$ thus yields an alternative formula for the knot longitude, which will be crucial later in generalising the Eisermann invariant. 
\end{Remark}

\subsubsection{One example of Eisermann invariants}
\begin{figure}
\centerline{\relabelbox 
\epsfysize 4.3cm 
\epsfbox{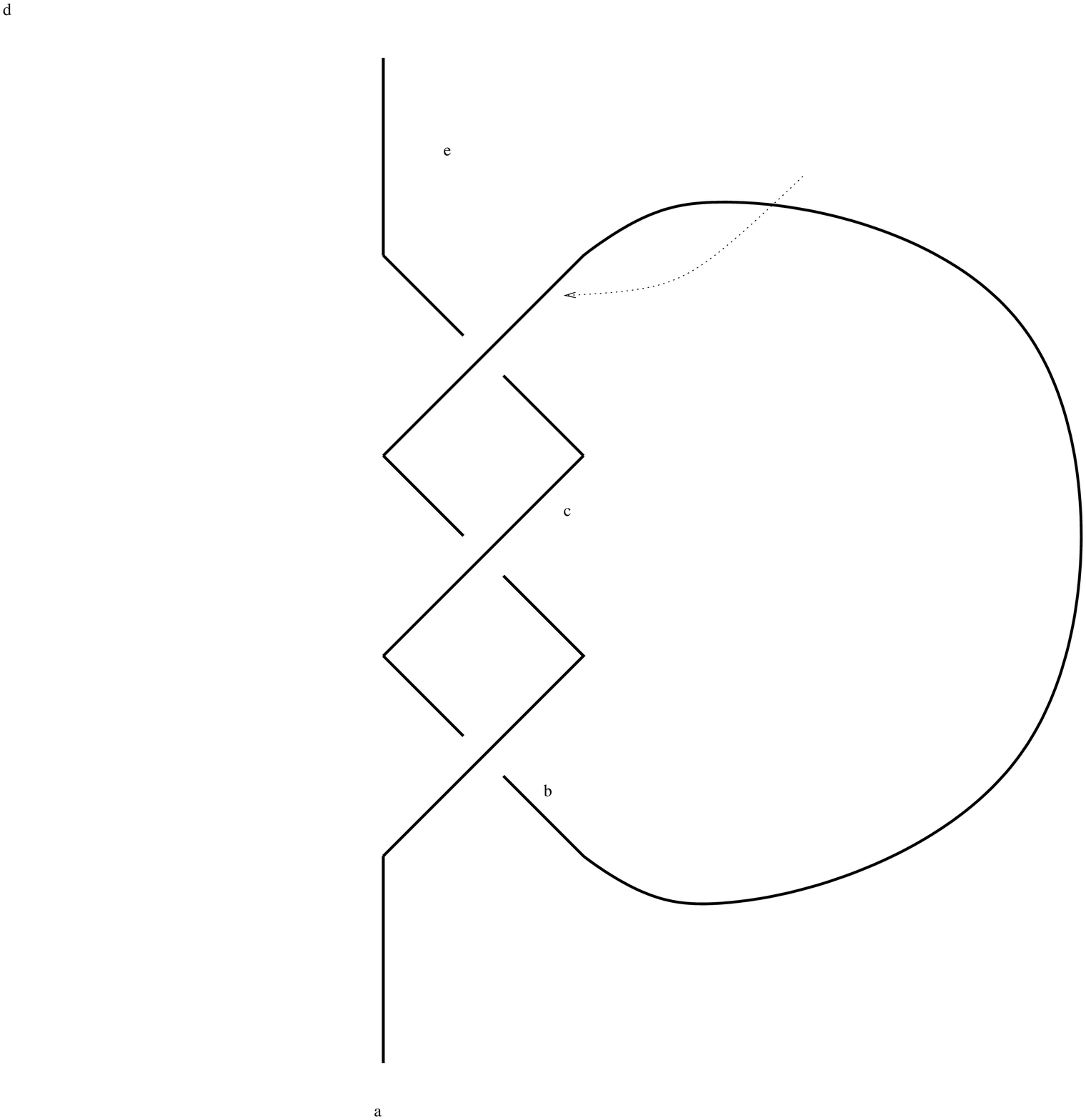}
\relabel{b}{${h}$}
\relabel{a}{${g}$}
\relabel{c}{${xhg^{-1}x^{-1}g}$}
\relabel{e}{${x^{2}gh^{-1}x^{-1}hg^{-1}x^{-1}g}$}
\relabel{d}{$x^{3}hg^{-1}x^{-1}gh^{-1}x^{-1}hg^{-1} x^{-1} g $}
\endrelabelbox}
\caption{Calculation of the Eisermann polynomial for the positive trefoil knot $K_+$\label{TP}}
\end{figure}

\begin{figure}
\centerline{\relabelbox 
\epsfysize 4.3cm 
\epsfbox{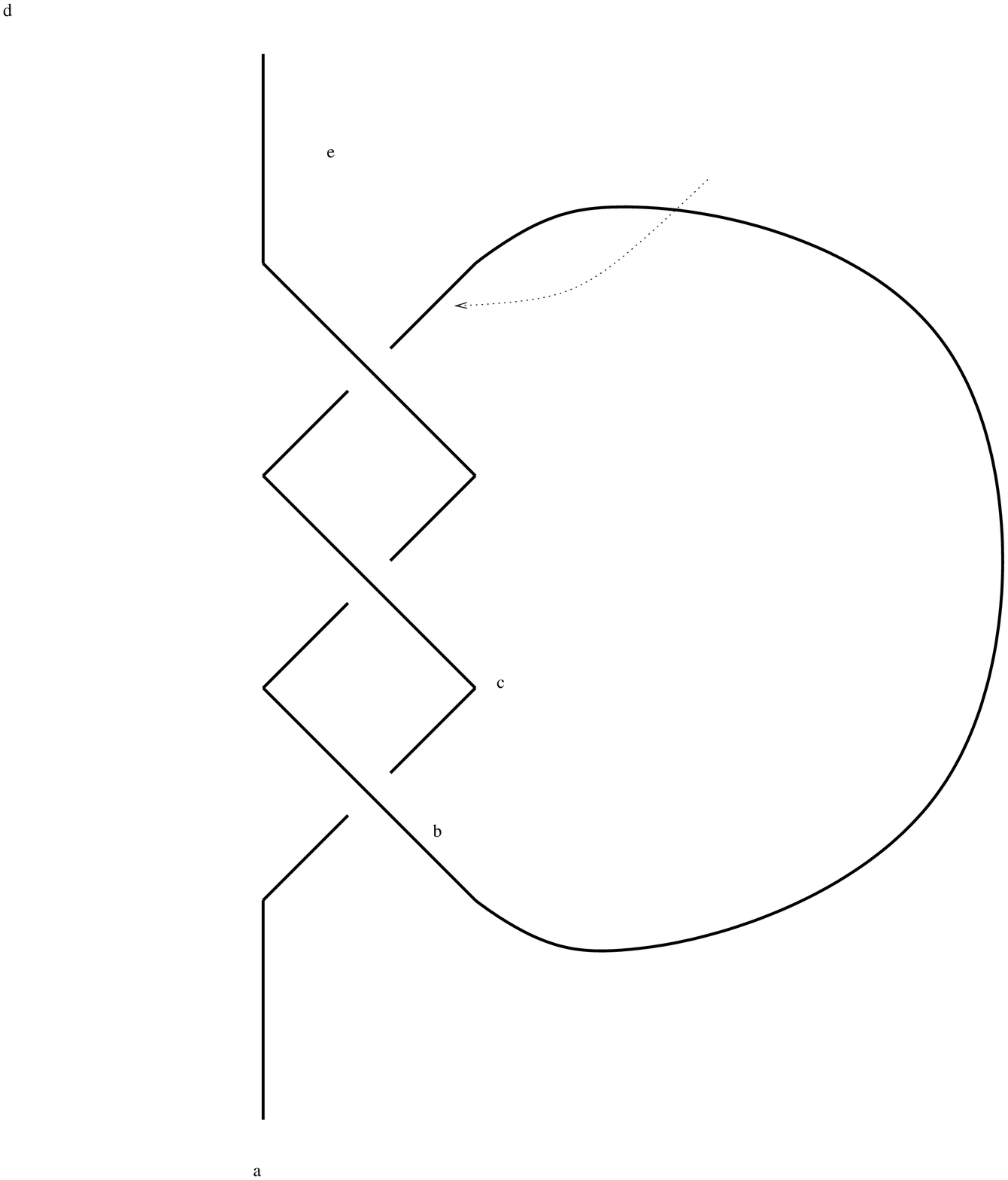}
\relabel{a}{${h}$}
\relabel{b}{${g}$}
\relabel{c}{${x^{-1}hg^{-1}xg}$}
\relabel{d}{${x^{-2}gh^{-1}xhg^{-1}xg}$}
\relabel{e}{$x^{-3}hg^{-1}xgh^{-1}xhg^{-1} x g $}
\endrelabelbox}
\caption{Calculation of the Eisermann polynomial for the negative trefoil knot $K_-$\label{TM}}
\end{figure}

Let $G$ be a group with a base point $x$.
The explicit calculation of the invariant $\langle a |I_{\bar{\Phi}^x}(K_+)| b \rangle $ and $\langle a |I_{\bar{\Phi}^x}(K_-)| b \rangle$
for the trefoil knot $K_+$ and its mirror image $K_-$ (the positive and negative trefoils), converted to string knots,  appears in figures \ref{TP} and \ref{TM}. In particular, given $a \in G$, we have:
\begin{multline*}
\langle a |I_{\bar{\Phi}^x}(K_+)| g \rangle=\\\# \big\{h \in G\colon   x^{3}hg^{-1}
x^{-1}gh^{-1}x^{-1}hg^{-1} x^{-1} g=a \textrm{ and }  x^{2}gh^{-1}x^{-1}hg^{-1}x^{-1}g=h \big\}, \, g \in G; 
\end{multline*}
\begin{equation*}
\langle a |I_{\bar{\Phi}^x}(K_-)| h \rangle=\# \big\{g \in G\colon   x^{-3}hg^{-1}xgh^{-1}xhg^{-1} x g=g \textrm{ and }  x^{-2}gh^{-1}xhg^{-1}xg=a \big\}, \, h \in G. 
\end{equation*}

Consider from now on $G=S_5$. We refer to   table \ref{T}, displaying the values of $\langle 1 |I_{\bar{\Phi}^x}(K_-)| a \rangle$ and  $\langle 1 |I_{\bar{\Phi}^x}(K_+)| a \rangle$ for some choices of $x \in S_5$, representing all possible conjugacy classes in $S_5$. In each second and third rows of table \ref{T}, we have $$\sum_{h \in S_5} \langle 1 |I_{\bar{\Phi}^x}(K_-)| h\rangle \textrm{ and }  \sum_{g \in S_5} \langle 1 |I_{\bar{\Phi}^x}(K_+)| g\rangle,$$ both elements of the group algebra of $S_5$.
\begin{table}
\begin{tabular}{|l|l|l|l|l|l|l|l}
 $x$  & $\id_{S_5}$ &$(12)$         & $(12)(34)$       & ($123)$         & $(123)(45)$ & $(1234)$  & $(12345)$ \\ \hline
 $K_-$ & $\id_{S_5}$& $7\id_{S_5}$& $5\id_{S_5}$   & $7\id_{S_5}$  &$\id_{S_5}$& $\id_{S_5}+4(13)(24)$   &$\id_{S_5}+5(12345)$  \\                      \hline 
 $K_+$ & $\id_{S_5}$& $7\id_{S_5}$& $5\id_{S_5}$   & $7\id_{S_5}$  &$\id_{S_5}$ & $\id_{S_5}+4(13)(24) $ & 
$\id_{S_5}+5(15432)$ 
\end{tabular}
\vskip 0.2 cm
\caption{\label{T} Value of the Eisermann invariant for the negative and positive trefoils for $G=S_5$.}
\end{table}
We can see that the invariant  $L \mapsto \langle 1 |I_{\bar{\Phi}^x}(L)| b \rangle$, where we identify a knot with its associated string knot, distinguishes the trefoils, for $x=(12345)$. This example is due to Eisermann \cite{E2}.

\subsection{2-crossed modules and braided crossed modules}
Two rich classes of examples of Reidemeister pairs  can be derived from 2-crossed modules. This will be the subject of  the end of this subsection and also \ref{LEI}. 

\subsubsection{Definition of 2-crossed modules}
A {pre-crossed module} $(\d\colon E \to G,\tr)$, is given by a group morphism $\d\colon E \to G$, together with a left action $\tr$ of $G$  on $E$ by automorphisms,  such that:  $$\d(g \tr e)=g\d(e)g^{-1}, \textrm{ for each } g \in G \textrm{ and each } e \in E.$$
This is the first Peiffer relation, appearing in the definition of a crossed module.

Let $(\d\colon E \to G, \tr)$ be a pre-crossed module.
Given $e,f  \in E$,  their {Peiffer commutator} is, by definition:
$$\langle  e,f\rangle   =\big(efe^{-1}\big)\,\, \big(\d(e) \tr f^{-1}\big).$$
Thus a  pre-crossed module  is a crossed module if and only if all of its Peiffer commutators vanish.

\begin{Definition}[(2-crossed module)]\label{2cmlg}
A 2-crossed module \cite{Co} is given by {a chain complex} of  groups:
$$L \ra{\de} E\ra{\d} G,$$
together with left actions $\tr$ by automorphisms of $G$ on $L$ and  $E$ (and on $G$ by conjugation), and a $G$-equivariant function {$\left\{,\right \} \colon E \times E \to L$} (called the Peiffer lifting). {Here $G$-equivariance means:}
\begin{equation*}
a \tr \{e,f\}= \{a \tr e, a \tr f\}, \fo a \in G \an e,f \in E.
\end{equation*}
{These should satisfy:}
\begin{enumerate}
 \item $L \ra{\de} E\ra{\d} G$ is {a chain complex} of $G$-modules (in other words $\d$ and $\de$ are $G$-equivariant {and $\d \circ \de =1$}.)
\item $\de(\left\{e,f\right \})=\left <e,f\right >,$ for each $e,f \in E$. 
\item $[l,k]=\left\{\de(l),\de(k)\right \}, $ for each $l,k \in L$. Here $[l,k]=lkl^{-1}k^{-1}$.
\item $\left\{\de(l),e\right \}\left\{e, \de(l)\right \}=l(\d(e) \tr l^{-1})$, for each $e\in E$ and $l \in L$.
\item $\left\{ef,g\right \}=\left\{e,fgf^{-1}\right \}\d(e)\tr \left\{f,g\right \}$, for each $e,f,g \in E$.
\item $\left\{e,fg\right \}=\left \{e,f\right\}(\d(e) \tr f) \tr' \left \{e,g\right\}$, where {$e,f,g \in E$}.
\end{enumerate}

\end{Definition} 
In the final equation above we have put:
\begin{equation}\label{tprime}
e \tr'l \doteq l\left\{\de(l)^{-1},e\right\}, \wh l\in L \an e \in E.
\end{equation}
It is proven in \cite{Co,GFM} that $\tr'$ is a left action of $E$ on $L$ by automorphisms and that,
 together with the map $\de\colon L \to M$, the action $\tr'$ by automorphisms defines a crossed module. In particular we have the second Peiffer identity
\begin{equation}\label{peif2de}
\de(l)\tr' m = l m l^{-1}.
\end{equation}
In \cite{Co,FMPi} it is shown that the following identity holds:
\begin{equation}\label{ef,g}
\{ef,g\} = (e \tr' \{f,g\})\{e,\d(f)\tr g\}.
\end{equation}
Applying (v),(vi) of Definition \ref{2cmlg} to $\{1_E1_E,e\}$ and $\{e,1_E1_E\}$ we obtain
\begin{equation}\label{1E}
\{1_E,e\}=\{ e, 1_E\} =1_L.
\end{equation}
Using condition (iv) of Definition \ref{2cmlg} and \eqref{tprime}, we also obtain:
\begin{equation}\label{abc}
\d(e) \tr l=(e\tr'{l})\,\, \{e,\de(l)^{-1}\}.
\end{equation}
\subsubsection{A technical lemma on 2-crossed modules}
In the following lemma we prove two algebraic equations which correspond to the Reidemeister 2 and 3 conditions.

\begin{Lemma}\label{rnn}In a 2-crossed module  $(L\ra{\de} E \ra{\d} G, \tr, \{,\})$  we have, for each $x,y,z \in E$:
\begin{equation}\label{invformula}
\begin{split}
x\tr' \{x^{-1},y\}\,\, \{x, \d(x)^{-1} \tr y \}&=1,\\
\{x,y\}\,\, (\d(x) \tr y)\tr' \{x,y^{-1}\}=1.
\end{split}
\end{equation}
\begin{equation}\label{mir}
\{x,y\} \,\, (\d(x) \tr y) \tr' \{x,z\}\,\,\{\d(x) \tr y , \d(x) \tr z\}=x \tr' \{y,z\}\,\, \{x, \d(y) \tr z\}\,\,(\d(xy) \tr z)\tr'\{x,y\}
\end{equation}
\end{Lemma}

\begin{Proof}
The first equation of \eqref{invformula} follows from \eqref{ef,g}, substituting $e=x,\, f=x^{-1},\, g=y$ and using \eqref{1E}. The second equation is analogous. To show \eqref{mir},
we apply condition (vi) of Definition \ref{2cmlg} to the underlined factors in:
$$\underline{\{x,y\} \,\, (\d(x) \tr y) \tr' \{x,z\}}\,\,\{\d(x) \tr y , \d(x) \tr z\}=x \tr' \{y,z\}\,\, \underline{\{x, \d(y) \tr z\}\,\,(\d(xy) \tr z)\tr'\{x,y\}}
$$
(with $e=x,\, f=y, \, g=z$ and $e=x,\, f=\d(y)\tr z, \, g=y$, respectively), to obtain:
$$\{x,yz\}\,\,\{\d(x) \tr y , \d(x) \tr z\}=\underline{x \tr' \{y,z\}} \,\,\{x,\d(y) \tr z\,\,y\}. $$
Using  \eqref{abc}, with $e=x, \, l = \{y,z\}$,  and condition (ii) of Definition \ref{2cmlg} on the underlined factor above, we obtain:
$$\{x,yz\}\,\,\{\d(x) \tr y , \d(x) \tr z\}=\d(x) \tr \{y,z\} \,\,\{x, \langle  y,z\rangle ^{-1}  \}^{-1} \,\,\{x,\d(y) \tr z\,\,y\} $$
or:
$$\{\d(x) \tr y , \d(x) \tr z\}=\underline{\{x,yz\}^{-1}\,\, \d(x) \tr \{y,z\} \,\,\{x, \langle  y,z\rangle ^{-1}  \}^{-1} }\,\,\{x,\d(y) \tr z\,\,y\}. $$
Applying \eqref{peif2de} to the underlined factor above, with $l= \d(x)\tr \{y,z\}^{-1},\, m=\{x,yz\}^{-1}$, together with the $G$-equivariance of $\de$, gives:
\begin{eqnarray*}
\lefteqn{\{\d(x) \tr y , \d(x) \tr z\}  =\d(x) \tr \{y,z\}\, . }\\  & & \quad \quad \quad \quad\,\,(\d(xy)\tr z\,\, \d(x) \tr (yz^{-1}y^{-1}))\tr'\{x,yz\}^{-1}\,\,\{x, \langle  y,z\rangle^{-1}\}^{-1} \,\,\{x,\d(y) \tr z\,\,y\} 
\end{eqnarray*}
or:
\begin{eqnarray*}
\lefteqn{\{\d(x) \tr y , \d(x) \tr z\}  =\d(x) \tr \{y,z\}\, . }\\  & & \quad \quad \quad \quad \underline{\,\,(\d(x)\tr \langle  y,z\rangle^{-1}))\tr'\{x,yz\}^{-1}\,\,\{x, \langle  y,z\rangle^{-1}\}^{-1}} \,\,\{x,\d(y) \tr z\,\,y\}  .
\end{eqnarray*}
Finally we apply condition (vi) of Definition \ref{2cmlg} to the underlined factor above, with $e=x,\, f= \langle  y,z\rangle^{-1}, \, g=yz$, taking the inverse equation and noting that $fg= \d(y) \tr z\,\,y$, to obtain:
$$\{\d(x) \tr y , \d(x) \tr z\}=\d(x) \tr \{y,z\} \,\,\{x,\d(y)\tr z\,\,y \}^{-1} \,\,\{x,\d(y) \tr z\,\,y\} = \d(x) \tr \{y,z\} $$
which holds because of the $G$-equivariance of {the Peiffer lifting} $\left\{,\right \}$.
\end{Proof}

\subsubsection{Braided crossed modules}

Later on we will need a generalisation of equation \eqref{mir}, for the case of braided crossed modules; see the example on page 55 of \cite{BG}.
\begin{Definition}[(Braided crossed module)]
A braided  crossed module is a 2-crossed module $(L\ra{\de} E \ra{\d} G, \tr, \{,\})$, where $G=\{1\}$, the trivial group. 
\end{Definition}
Note that in a braided crossed module we have, for each $e,f\in E$: 
\begin{equation}\label{pcbcm}
\de(\{e,f\})=[e,f]=efe^{-1}f^{-1}. 
\end{equation}
Also, it is proven in \cite{Co} and follows from equation (13) of \cite{GFM} that:
\begin{equation}\label{actform}
g \tr' \{e,f\}=\{geg^{-1},gfg^{-1}\}, \textrm{ for each } e,f,g \in E. 
\end{equation}

If $(L\ra{\de} E \ra{\d} G, \tr, \{,\})$ is a braided crossed module, then the underlying monoidal category of the categorical group $\C(\de \colon L \to E,\tr')$ is braided \cite{JS}, the braiding $gh=g \tn h\to h \tn g =hg$  being determined by the Peiffer pairing $\{g,h\}^{-1}$; {see page 56 of \cite{BG}. The framed knot invariants \cite{K,CP} associated to this ${\cal R}$-matrix and an object of $\C(\de \colon L \to E,\tr')$, in other words an element of $E$, are however trivial, since it is easy to see that $E$  acts trivially on $\ker(\de) \subset L$, through $a \tr' l=l\{\de(l),a\}$. Nevertheless we can modify this ${\cal R}$-matrix to an ${\cal R}$-matrix in a linearisation of $\C(\de \colon L \to E,\tr')$, which yields  non-trivial knot invariants; see Theorem \ref{el2}. Let us now prove some equations leading to this.}

Considering a braided crossed module, let us look at  equation \eqref{mir}, with $x=Mx^{-1}$, $y=Lx^{-1}$ and $z=Nx^{-1}$ (where $M,N,L,x \in E$), which explicitly is:
\begin{multline*}
\{Mx^{-1},Lx^{-1}\}\,\,L \tr' \{x^{-1} M,x^{-1}N\}\,\,\{Lx^{-1},Nx^{-1}\}=\\
M \tr' \{x^{-1}L,x^{-1}N\}\,\,\{Mx^{-1},Nx^{-1}\}\,\,N\tr'\{x^{-1}M,x^{-1}L\}. 
\end{multline*}
By using \eqref{actform}, \eqref{peif2de} and \eqref{pcbcm} this can be rearranged as:
\begin{multline}\label{squash}
\{Mx^{-1},Lx^{-1}\}\,\,\{Lx^{-1},Nx^{-1}\}\,\\\,\{Nx^{-1}LN^{-1}Mx^{-1}NL^{-1}xN^{-1},Nx^{-1}Lx^{-1}NL^{-1}xN^{-1}\}\\
= \{Mx^{-1},Nx^{-1}\} \,\,\\ \{ Nx^{-1}MN^{-1}Lx^{-1}NM^{-1}xN^{-1} ,   Nx^{-1}Mx^{-1}NM^{-1}xN^{-1}  \}\,\,(Nx^{-1})\tr'\{Mx^{-1},Lx^{-1}\}. 
\end{multline}
\begin{Lemma}In a braided crossed module, for $M,L,N,x \in E$ we have:
\begin{multline}\label{r3sharp}
 \{Mx^{-1},Lx^{-1}\}\,\,L\tr'\{x^{-1}ML^{-1}xLx^{-1},Nx^{-1}\}\,\,\{Lx^{-1},Nx^{-1}\}\\
=M\tr'\{Lx^{-1},Nx^{-1}\}\,\,\{Mx^{-1},Nx^{-1}\}\,\,N \tr' \{x^{-1}MN^{-1}xNx^{-1},x^{-1}LN^{-1}xNx^{-1}\}.
\end{multline}
\end{Lemma}
\begin{Proof}
 This is the same as, by using \eqref{actform}, {\eqref{peif2de} and \eqref{pcbcm}:}
\begin{multline*}
 \{Mx^{-1},Lx^{-1}\}\{Lx^{-1},Nx^{-1}\}\\ \{Nx^{-1}LN^{-1}ML^{-1}xLx^{-1}x^{-1}NL^{-1}xN^{-1}, Nx^{-1}LN^{-1}xNx^{-1}x^{-1}NL^{-1}xN^{-1}   \}\\
=\{Mx^{-1},Nx^{-1}\} \,\, \{Nx^{-1}MN^{-1}xLx^{-2}NM^{-1}xN^{-1},Nx^{-1}MN^{-1}xNx^{-2}NM^{-1}xN^{-1} \}\\ \,\,N \tr'\{x^{-1}MN^{-1}xNx^{-1},x^{-1}LN^{-1}xN x^{-1}\},
\end{multline*}
or:
\begin{multline*}
 \{Mx^{-1},Lx^{-1}\}\{Lx^{-1},Nx^{-1}\}\\  \{Nx^{-1}LN^{-1}M\,\,[L^{-1},x]\,\, x^{-1} NL^{-1}xN	^{-1}, Nx^{-1}L\,\,[N^{-1},x] \,\,x^{-1}NL^{-1}xN^{-1}   \}\\
=\{Mx^{-1},Nx^{-1}\} \,\, \{Nx^{-1}MN^{-1}L\,\,[L^{-1}, x]  \,\,  x^{-1}NM^{-1}xN^{-1},Nx^{-1}M \,\,[N^{-1},x]\,\,x^{-1}NM^{-1}xN^{-1} \}\\ N \tr' \{x^{-1}M\,\,[N^{-1},x],x^{-1} L\,\,[N^{-1},x]\}.
\end{multline*}
This follows, by \eqref{pcbcm}, from the more general (where $A,B\in L$):
\begin{multline}\label{toprove}
 \{Mx^{-1},Lx^{-1}\}\{Lx^{-1},Nx^{-1}\} \\ \{Nx^{-1}LN^{-1}M\,\,\de(A)\,\, x^{-1} NL^{-1}xN	^{-1}, Nx^{-1}L\,\,\de(B) \,\,x^{-1}NL^{-1}xN^{-1}   \}\\
=\{Mx^{-1},Nx^{-1}\} \,\, \{Nx^{-1}MN^{-1}L\,\,\de(A)  \,\,  x^{-1}NM^{-1}xN^{-1},Nx^{-1}M \,\,\de(B)\,\,x^{-1}NM^{-1}xN^{-1} \}\\ \{Nx^{-1}MN^{-1}\,\,\de(N \tr' B),Nx^{-1}LN^{-1}\,\,\de(N \tr' B)\}.
\end{multline}
By using \eqref{ef,g} and the final condition of the definition of a 2-crossed module,  together with \eqref{peif2de}, we can see that, in a braided crossed module,  where $a,b \in E$ and $k,l \in L$:
\begin{align*}
\{a \de(k),b \de(l)\}& =a \tr' k\,\, (ab) \tr' (lk^{-1}l^{-1})\,\,\{a,b\}\,\, (ba)\tr' l\,\,b \tr' l^{-1} \\
&=\{a,b\}\,\,(bab^{-1}) \tr' k\,\, (ba) \tr' (lk^{-1})\,\,b \tr' l^{-1}\\
&=a \tr' k\,\, (ab) \tr' (lk^{-1})\,\,(aba^{-1}) \tr' l^{-1}\,\,\{a,b\}.
\end{align*}
When $k=l$ this means:
\begin{align}
\{a \de(k),b \de(k)\}&=a \tr' k\,\,\{a,b\}\,\,b \tr' k^{-1}\label{eqa}
\\&=\{a,b\}\,\,(bab^{-1}) \tr' k\,\,b \tr' k^{-1}.
\end{align}

Also (where $a,a',b,b' \in E$ and $k,l \in L$):
\begin{equation}\label{eqb}
\begin{split}
\{a \de(k)a'&,b \de(l)b'\}=\left\{(aa')\de(a'^{-1} \tr' k),(bb')\de(b'^{-1} \tr' l)\right\}\\
&=\{aa',bb'\}\,\,(bb'aa'b'^{-1}b^{-1}a'^{-1}) \tr' k\,\, (bb'aa') \tr' (b'^{-1}\tr'l\,\, a'^{-1}\tr' k^{-1})\,\,b \tr' l^{-1}.
\end{split}
\end{equation}
Applying equations \eqref{eqa} and \eqref{eqb} to \eqref{toprove}, together with the 2-crossed module relation $k l=l\,\, \de(l)^{-1} \tr' k$, for each $k,l \in L$ makes \eqref{toprove} follow plainly from \eqref{squash}. 
\end{Proof}

\subsubsection{Peiffer liftings define Reidemeister pairs}

Let $(L \ra{\de} E \ra{\d} G, \tr, \{,\})$ be a 2-crossed module. It is easy to see that the operations: 
$$x \tr y \doteq \d(x) \tr y \textrm{ and } y \tl x\doteq  \d(x)^{-1}\tr y ,$$
where $x,y \in E$,
define a rack structure in $E$. Indeed, for example:
\begin{align*}
x \tr (y \tr z)=&\d(x) \tr (\d(y) \tr z)=\d(xy) \tr z=\d(xyx^{-1}) \tr (\d(x) \tr z) \\&=\d(\d(x) \tr y) \tr( \d(x) \tr z)\\
&=(x \tr y) \tr (x \tr z).
\end{align*}
Any finite 2-crossed module defines, as long as $E$ is finite, an invariant of framed knots by using the rack structure in $E$,  Theorem \ref{rli}. The Peiffer lifting permits us to refine this rack invariant with an element of $\ker(\de) \subset L$. 
\begin{Theorem} Suppose that $E$ is finite. 
Consider the crossed module $(\de\colon L \to E, \tr')$.
Put 
\begin{equation}
 \begin{split}
  \psi(A,B)&=\{A,B\},\\
  \phi(A,B)&=A \tr'\{A^{-1},B\}.
 \end{split}
\end{equation}
 Then $\Phi=(\psi,\phi)$ is a framed Reidemeister pair.
\end{Theorem}
\begin{Proof}
 Equation $R_2$ of Definition \ref{rp} is in this case:
$$\phi(A,B)\,\,\psi(A,\d(A)^{-1}\tr B)=1,$$
or
$$A \tr' \{A^{-1},B\}  \,\, \{A,\d(A)^{-1} \tr B\}=1, $$
which follows at once from Lemma \ref{rnn}.
Equation \eqref{R3p} is exactly the second equation of Lemma \ref{rnn} with $x=X$, $y=Y$ and $z=Z$.

The equation $\de(\phi(A,Z)) A=Z$ means exactly that $Z=\d(A) \tr A$. By the Nelson Lemma \ref{Nelson} $A \mapsto \d(A) \tr A$ is bijective, thus for each $Z$ the equation $\d(\phi(A,Z)) A=Z$ has a unique solution $f(Z)$. Moreover $g(A)=\de(\psi(A,A))^{-1}A=\d(A) \tr A$, thus $g$ is the inverse of $f$.

\end{Proof}

 Let $G$ be a finite group. Let $G^{\ab}$ be its abelianisation. Consider the 2-crossed module: $$(G^\ab \otimes_{\Z} G^\ab \ra{\de} G \ra{\id} G),$$ where the action of $G$ on $G$ is by conjugation and on $G^\ab \tn_\Z G^\ab$ is trivial. Also $\de(G \tn_{\Z} G)=\{1_G\}$.
 The Peiffer pairing is $$\{a,b\}=a^\ab \tn b^\ab.  $$Here $x\in  G \mapsto x^\ab \in G^\ab$ is the abelianisation map. Given a knot $K$, the associated knot invariant essentially counts the number of representations of the knot group into $G$, separating then according to the image of an arc into $G^{\ab}$, squared. Recall that the abelianisation of a knot group $\pi_1(C_K)$ is always $\Z$ and that all elements of  $\pi_1(C_K)$ associated to the oriented arcs of {a diagram of} $K$ project to the same element in the abelianisation. 

We thus have:
\begin{Theorem}
 Let $G$ be a finite group. Let $\Phi$ be the framed Reidemeister pair given by the 2-crossed module $(G^\ab \otimes G^\ab \ra{\de} G \ra{\id} G)$. Let $K$ be a framed link. Then:
$$I_{\Phi}(K)=\sum_{f\colon \pi_1(C_K) \to G} \big (f(m)\tn f(m)\big)^t \in \Z[G^{\ab} \otimes_\Z G^{\ab}], $$
where $m \in \pi_1(C_K)$ is any meridian of $K$ and $t$ is the framing of $K$.
\end{Theorem}

\subsection{Lifting the Eisermann invariant}\label{LEI}

Recall the construction of the Eisermann invariant  in our framework (subsection \ref{ep}). The quandle underlying the Eisermann invariant corresponds in our approach to a certain unframed Reidemeister pair given in equation 
(\ref{eisphipsi}) by:
\begin{eqnarray*}
\phi^x(L,M) &= &  [Mx^{-1},Lx^{-1} ],\\
\psi^x(L,M)&= & [L,M][ML^{-1},x]=[xML^{-1}x^{-1}Lx^{-1},Lx^{-1}]^{-1}.
\end{eqnarray*}
Since the expressions for the functions {$\phi^x$ and $\psi^x$} are given in terms of commutators, it is natural, given \eqref{pcbcm}, to seek generalizations of them in the context of braided crossed modules.
\subsubsection{The unframed case}\label{tuc}
\begin{Definition}[(Unframed Eisermann lifting) ] Let $G$ be a group and $x \in G$.
An  unframed Eisermann lifting is given by a crossed module $(\d: E\rightarrow G, \tr)$, and an unframed Reidemeister pair $\Phi^x=(\phi^x, \psi^x)$,  where $\phi^x, \, \psi^x: G \times G \rightarrow E$, such that: 
\begin{eqnarray*}
\d(\phi^x(L,M)) &= & [Mx^{-1},Lx^{-1} ],\\
\d(\psi^x(L,M)) &= & [L,M][ML^{-1},x].
\end{eqnarray*}
\end{Definition}
The colourings of the arcs of any tangle diagram will correspond to those given by the Eisermann quandle, but there may be additional information contained in the assignments of elements of $E$ to the crossings of the diagram, i.e. an Eisermann lifting is a refinement of the Eisermann invariant.

In this subsection we will provide examples of unframed Eisermann liftings coming from braided crossed modules and from central extensions of groups.
\begin{Example}\label{2xmodE}
Given a braided crossed module 
$$
( E \ra{\de} G \to 1, \{ \, , \, \})
$$
such that $\de$ is surjective, and $x \in G$, define $\Phi^x =(\psi^x, \phi^x)$ by:
$$
\phi^x(L,M) = \{ Mx^{-1}, Lx^{-1} \}, \quad \quad \psi^x(L,M) = \{L,M\} \{ ML^{-1}, x\}.
$$
\begin{Theorem}\label{el1}
{The pair $\Phi^x$} is an unframed Reidemeister pair for the crossed module $(\de\colon L \to G, \tr')$, obtained from the 2-crossed module $( L \ra{\de} G \ra{\d} 1, \{ \, , \, \})$, and constitutes an Eisermann lifting for this crossed module. 
\label{lift2Xm}
\end{Theorem}
\begin{Proof}
Given $L,M \in G$, then $\de (\{L,M\}) = [L,M]$ since the Peiffer commutator is the ordinary commutator in this case. Since $\de$ is surjective we can find $l\in L$ such that $L=\de(l)$, and likewise $M=\de(m)$ and $x=\de(y)$. Thus, using (3) of Definition \ref{2cmlg}, the Reidemeister 2 condition becomes:
$$
[my^{-1},ly^{-1}]\,\,[l, y^{-1}ml^{-1}yl]\,\,[y^{-1} ml^{-1}y,y]=1.
$$ 
This is an algebraic identity which was shown to hold in the proof of Theorem \ref{eisermann}. An analogous argument shows that the Reidemeister 3 equation is satisfied, since $\de(l)\tr' m = lml^{-1}$, so that we can use the algebraic identity for Reidemeister 3 from the same proof. For the Reidemeister 1 move this follows from 
$$\psi^x(M,M)=\{Mx^{-1},Mx^{-1}\}=\{\de(my^{-1}),\de(my^{-1})\}=[my^{-1},my^{-1}]=1_E .$$
\end{Proof}
\begin{Remark}
 By the proof of the previous theorem, we can see that an alternative expression for $\psi^x$ is:
$$\psi^x(L,M)=\{xML^{-1}x^{-1}Lx^{-1},Lx^{-1}\}^{-1}.$$
Note that for each $L,M,x \in E$ we have $\{xML^{-1}x^{-1}Lx^{-1},Lx^{-1}\}^{-1}=\{L,M\} \{ ML^{-1}, x\}$, since these are identities between usual commutators. 
\end{Remark}

\end{Example}

Next we construct Eisermann liftings using central extensions of groups. We first define a  2-crossed module given any central extension $\{0\}\to A \to E \to G \to \{1\}$.

Let $G$ be a group. Let $\d \colon E \to G$ be a surjective group morphism, such that the kernel $A$ of {$\d$ is central in $E$}. Clearly: 
\begin{Lemma}
Given $a,b \in E$, if $\d(a)=\d(b)$ then 
\begin{itemize}
\item $aea^{-1}=beb^{-1}$, for each $e \in E$. 
\item $[a,e]=[b,e]$, for each $e \in E$.
\item $[e,a]=[e,b]$, for each $e \in E$.
\end{itemize}
\end{Lemma}

Choose an arbitrary section $s \colon G \to E$ of $\d$. Therefore $s(ab)=s(a)s(b)\lambda(a,b)$ where $\lambda(a,b)$ in the centre of $E$. Clearly:

\begin{Lemma}\label{liftextension}
 The map $$(g,e) \in G \times E \mapsto g \tr e=s(g) \,\, e\,\, s(g)^{-1}\in E$$ is a left action of $G$ on $E$ by automorphisms, and with this action $\d\colon E \to G$ is a crossed module. Moreover, the action $\tr$ does not depend on the chosen section $s$.
\end{Lemma}
Note that $s$ need not be a group morphism. However:
\begin{Lemma}
 For each $g,h \in G$ and $a \in E$. then:
\begin{itemize}
\item $[s(gh),a]=[s(g)s(h),a] $,
\item $[s(g)^{-1},a]=[s(g^{-1}),a]$,
\item $[s([g,h]),a]=[[s(g),s(h)],a] $,
\end{itemize}
and the equalities obtained by swapping the two arguments in the commutators. Furthermore, if $a,b \in E$: 
$$[s(\d(a)),s(\d(b)]=[a,b].$$
\end{Lemma}
\begin{Proof}
 Certainly $s(gh)=\lambda(g,h)\, s(g)\, s(h)$, for each $g,h \in G$, where $\lambda(g,h) \in  A$. Thus 
$$[s(gh),a]=[s(g)s(h)\lambda(g,h),a]=[s(g)\, s(h),a].$$
The other equalities follow analogously.
\end{Proof}

{Given a section $s\colon G \to E$ of $\partial\colon E \to G$,} define, for each $g,h \in G$:
\begin{equation}\label{dpl}
 \{g,h\}=[s(g),s(h)]. 
\end{equation}
{(This does not depend on the chosen section since $\ker(\partial)$ is central.)}

\begin{Lemma}
For each $a,b \in E$: $$[a,b]=\{\d(a),\d(b)\} .$$ 
\label{bracketlift}
\end{Lemma}
\begin{Proof}
Given $a,b \in E$  $$\{\d(a),\d(b)\}=[s(\d(a)),s(\d(b)]=[a,b]. $$ 
 
\end{Proof}

\begin{Lemma}With this Peiffer lifting we have a 2-crossed module: $$( E \ra{\d} G \to \{1\},\{,\}) .$$
Moreover $\partial$ is surjective. 

\end{Lemma}
\begin{Proof}
 Condition (i) is immediate. 
 Condition (ii) is proven as follows:
$$\d(\{g,h\})=\d([s(g),s(h)])=[\d(s(g)),\d(s(h))]= [g,h].$$

Condition (iii) is Lemma \ref{bracketlift}.

Condition (iv). Given $a\in E$ and $g \in G$:
\begin{align*}
 \{\d(a),g\}\,\{g,\d(a)\}=[s(\d(a)),s(g)]\, [s(g),s(\d(a))]=1.
\end{align*}

Condition (v). Given $e,f,g \in G$ then:
\begin{align*}
 \{ef,g\}=[s(ef),s(g)]=[s(e)\,s(f),s(g)]
\end{align*}
and:
\begin{align*}
 \{e,fgf^{-1}\}\,\{f,g\}=[s(e),s(fgf^{-1})]\,[s(f),s(g)]=[s(e),s(f)\,s(g)\,s(f)^{-1}]\,[s(f),s(g)]=[s(e)\,s(f),s(g)].
\end{align*}

Condition (vi). Given $k \in E$ and $g \in G$ then: 
$$g \tr' k=k\{\d(k)^{-1},g\}=k[s(\d(k^{-1})),s(g)]=k [k^{-1},s(g)]=s(g)\, k\,s(g)^{-1} = g\tr k,$$
therefore:
\begin{align*}
\{e,f\} f \tr' \{e,g\}=[s(e),s(f)]\, s(f)\, [s(e),s(g)] s(f)^{-1}=[s(e),s(f)\,s(g)]=[s(e),s(fg)]=\{e,fg\}.
\end{align*}

\end{Proof}

Therefore:
\begin{Corollary}\label{liftextension2} Let $G$ be a group and $x \in G$. Let $\d\colon E \to G$ be a surjective group morphism such that the kernel $A$ of $\partial$ is central in $E$. The pair $\Phi^x =(\phi^x, \psi^x)$, given by
$$
\phi^x(g,h)=\{hx^{-1}, gx^{-1}\}, \quad \psi^x(g,h) = \{g,h\} \{hg^{-1},x\},
$$
{in terms of the Peiffer lifting \eqref{dpl},} is an unframed Eisermann lifting for the crossed module $(\d \colon E \to G, \tr)$, where $\tr$ is the left {action of $G$ on $E$} of Lemma \ref{liftextension}.
\end{Corollary}

\begin{Remark}
A class of central extensions for the preceding construction is given by the non-abelian wedge square of a finite group $G$, namely the map $\de: G\wedge G \rightarrow G'$, where $\d(g\wedge h) = [g,h]$ for each $g,h \in G$. The topological significance of these central extensions was emphasized by Brown and Loday \cite{BrL}. We will use these notions later when interpreting Eisermann liftings.
\end{Remark}

\subsubsection{A non-trivial example of {an} unframed  Eisermann lifting}\label{ufel}
In the context of {Corollary \ref{liftextension2},} let us find a lifting of the Eisermann invariant for the case of $G=S_5$, for which we gave detailed calculations in subsection \ref{ep}. It is well known that $S_5$ is isomorphic to $\PGL(2,5)$, the group of invertible two-by-two matrices in the field $\Z_5$, modulo the central subgroup $\Z_5^\ast$ of diagonal matrices which are multiples of the identity. We thus have a central extension:
$$\{0\} \to \Z_5^\ast \ra{i} \GL(2,5) \ra{p}  \PGL(2,5))\cong S_5\to \{1\}.$$ 
Here $\GL(2,5)$ is the group of invertible two-by-two matrices in the field $\Z_5$.

Let $K_+$ and $K_-$ be the right and left handed trefoils. In  table \ref{T2} we display $\langle a |I_{\Phi^x}(K_+)| 1 \rangle$ and  $\langle a |I_{\Phi^x}(K_-)| 1 \rangle$ for some choices of $x \in \PGL(2,5)\cong S_5$, representing all possible conjugacy classes in $S_5\cong \PGL(2,5)$, in the same order as in table \ref{T}. In the second and third rows of the table,  we have $$\sum_{s \,\in\, \PGL(2,5)} \langle s |I_{{\Phi}^x}(K_+)| 1\rangle \textrm{ and } \sum_{s \,\in \,\PGL(2,5)} \langle s |I_{{\Phi}^x}(K_-)| 1\rangle,$$
respectively, both elements of the group algebra of $\GL(2,5)$. If $A\in \GL(2,5)$, its projection to $\PGL(2,5)$ is denoted by $\widetilde{A}$. 

\begin{table}
\scriptsize
$$
{
\begin{array}{|l|l|l|l|l|l|l|l|}\hline
 x    &  \widetilde{\begin{pmatrix}  1 & 0\\ 0 &1 \end{pmatrix}}   &   \widetilde{ \begin{pmatrix}  1 & 3\\ 4 & 4  \end{pmatrix}}  &   \widetilde{\begin{pmatrix}  0 & 1\\ 1 & 0 \end{pmatrix}}       &   \widetilde{\begin{pmatrix}  4 & 1\\ 4 & 0  \end{pmatrix}}         &  \widetilde{\begin{pmatrix}  3 & 1\\ 4 &4 \end{pmatrix}  }    &   \widetilde{ \begin{pmatrix}  2 & 0\\ 0 & 1 \end{pmatrix}  }    &    \widetilde{\begin{pmatrix}  3 & 0\\ 3 & 3 \end{pmatrix} }      \\ \hline
 K_+  & \begin{pmatrix} 1 & 0\\ 0 &1 \end{pmatrix}    &   7\begin{pmatrix}  1 & 0\\ 0 &1 \end{pmatrix}   &  5 \begin{pmatrix}  1 & 0\\ 0 & 1 \end{pmatrix}     & \begin{pmatrix}  1 & 0\\ 0 & 1  \end{pmatrix}   + 6\begin{pmatrix}  4 & 0\\ 0 & 4  \end{pmatrix}           &\begin{pmatrix}  1 & 0\\ 0 & 1 \end{pmatrix} &  \begin{pmatrix}  1 & 0\\ 0 & 1 \end{pmatrix}  +4 \begin{pmatrix}  3 & 0\\ 0 & 2 \end{pmatrix}   &  \begin{pmatrix}  1 & 0\\ 0 & 1 \end{pmatrix}  +5 \begin{pmatrix}  4 & 0\\ 1 & 4 \end{pmatrix}     \\                      \hline 
 K_-  & \begin{pmatrix}  1 & 0\\ 0 &1 \end{pmatrix}   &   7\begin{pmatrix}  1 & 0\\ 0 &1 \end{pmatrix} &   5 \begin{pmatrix}  1 & 0\\ 0 & 1 \end{pmatrix}      &  \begin{pmatrix}  1 & 0\\ 0 & 1  \end{pmatrix}   + 6\begin{pmatrix}  4 & 0\\ 0 & 4  \end{pmatrix}  & \begin{pmatrix}  1 & 0\\ 0 & 1 \end{pmatrix}   &   \begin{pmatrix}  1 & 0\\ 0 & 1 \end{pmatrix}  +4 \begin{pmatrix}  2 & 0\\ 0 & 3 \end{pmatrix}   &   \begin{pmatrix}  1 & 0\\ 0 & 1 \end{pmatrix}  +5 \begin{pmatrix}  4 & 4\\ 4 & 0 \end{pmatrix}\\\hline    
\end{array}}
$$
\normalsize
\caption{\label{T2} Value of the lifted Eisermann invariant for the positive and negative trefoils in subsection \ref{ufel}}
\end{table}
\begin{table}\scriptsize
$$
\begin{array}{|l|l|l|l|l|l|l|l|}\hline
 x    &  \widetilde{\begin{pmatrix}  1 & 0\\ 0 &1 \end{pmatrix}}   &   \widetilde{ \begin{pmatrix}  1 & 3\\ 4 & 4  \end{pmatrix}}  &   \widetilde{\begin{pmatrix}  0 & 1\\ 1 & 0 \end{pmatrix}}       &   \widetilde{\begin{pmatrix}  4 & 1\\ 4 & 0  \end{pmatrix}}         &  \widetilde{\begin{pmatrix}  3 & 1\\ 4 &4 \end{pmatrix}  }    &   \widetilde{ \begin{pmatrix}  2 & 0\\ 0 & 1 \end{pmatrix}  }    &    \widetilde{\begin{pmatrix}  3 & 0\\ 3 & 3 \end{pmatrix} }      \\ \hline
 K_+  & \widetilde{\begin{pmatrix} 1 & 0\\ 0 &1 \end{pmatrix}}    &   7\widetilde{\begin{pmatrix}  1 & 0\\ 0 &1 \end{pmatrix}}   &  5 \widetilde{\begin{pmatrix}  1& 0\\ 0 & 1 \end{pmatrix} }    & 7\widetilde{\begin{pmatrix}  1 & 0\\ 0 & 1  \end{pmatrix}}           &\widetilde{\begin{pmatrix}  1 & 0\\ 0 & 1 \end{pmatrix}} &  \widetilde{\begin{pmatrix}  1 & 0\\ 0 & 1 \end{pmatrix}}  +4 \widetilde{\begin{pmatrix}  4 & 0\\ 0 & 1 \end{pmatrix}}   & \widetilde{ \begin{pmatrix}  1 & 0\\ 0 & 1 \end{pmatrix}}  +5 \widetilde{\begin{pmatrix}  3& 0\\ 3 & 3 \end{pmatrix} }    \\                      \hline 
 K_-  & \widetilde{\begin{pmatrix}  1 & 0\\ 0 &1 \end{pmatrix}}   &   7 \widetilde{\begin{pmatrix}  1 & 0\\ 0 &1 \end{pmatrix}} &   5 \widetilde{\begin{pmatrix}  1 & 0\\ 0 & 1 \end{pmatrix}}      &  7\widetilde{\begin{pmatrix}  1 & 0\\ 0 & 1  \end{pmatrix}}  & \widetilde{\begin{pmatrix}  1 & 0\\ 0 & 1 \end{pmatrix}}   &  \widetilde{ \begin{pmatrix}  1 & 0\\ 0 & 1 \end{pmatrix}}  +4 \widetilde{\begin{pmatrix}  4 & 0\\ 0 & 1 \end{pmatrix}}   & \widetilde{  \begin{pmatrix}  1 & 0\\ 0 & 1 \end{pmatrix}}  +5\widetilde{ \begin{pmatrix}  2& 0\\ 3 & 2 \end{pmatrix}}    \\ \hline
\end{array}
$$
\normalsize
\caption{\label{Ts} Value of the unlifted Eisermann invariant for the positive and negative trefoils in subsection \ref{ufel}}
\end{table}

Comparing with table \ref{Ts}, which shows the unlifted {Eisermann invariant $I_{\Phi^x_0}$ for $G={\rm PGL}(2,5)$,}  we can see that this lifting of the Eisermann invariant is strictly stronger than the Eisermann invariant itself. Specifically, looking at the penultimate column of tables \ref{T2} and \ref{Ts}, thus $x=\widetilde{ \begin{pmatrix}  2 & 0\\ 0 & 1 \end{pmatrix}  } $, we can see that the lifting distinguishes the trefoil from its mirror image. Namely, for the lifted one, noting that  $\widetilde{\begin{pmatrix} 3 &0\\0 &2 \end{pmatrix}}= \widetilde{\begin{pmatrix} 2 &0\\0 &3\end{pmatrix}}=\widetilde{\begin{pmatrix} 4 &0\\0 &1\end{pmatrix}}$, we have:
$$\left \langle \widetilde{\begin{pmatrix} 3 &0\\0 &2\end{pmatrix}} \Big | I_{{\Phi}^x}(K_+) \Big |  \widetilde{\begin{pmatrix} 1 & 0 \\ 0  & 1   \end{pmatrix} } \right \rangle=4\begin{pmatrix} 3 &0\\0 &2 \end{pmatrix}\neq 4\begin{pmatrix} 2 &0\\0 &3 \end{pmatrix} =\left \langle \widetilde{\begin{pmatrix} 3 &0\\0 &2\end{pmatrix}} \Big | I_{{\Phi}^x}(K_-) \Big |  \widetilde{\begin{pmatrix} 1 & 0 \\ 0  & 1   \end{pmatrix} } \right \rangle,$$                    
whereas for the unlifted Eisermann invariant $I_{\Phi^x_0}$:
$$\left \langle \widetilde{\begin{pmatrix} 3 &0\\0 &2\end{pmatrix}} \Big | I_{\Phi^x_0}(K_+) \Big |  \widetilde{\begin{pmatrix} 1 & 0 \\ 0  & 1   \end{pmatrix} } \right \rangle=4\widetilde{\begin{pmatrix} 3 &0\\0 &2 \end{pmatrix}}= \left \langle \widetilde{\begin{pmatrix} 3 &0\\0 &2\end{pmatrix}} \Big | I_{\Phi^x_0}(K_-) \Big |  \widetilde{\begin{pmatrix} 1 & 0 \\ 0  & 1   \end{pmatrix} } \right \rangle.$$                    
Table \ref{Ts} should be compared with table \ref{T}. 

\subsubsection{Framed Eisermann liftings from braided crossed modules}\label{fel}

\begin{Definition}[Framed Eisermann lifting ] Let $G$ be a group and $x \in G$.
A  framed Eisermann lifting is given by a crossed module $(\d: E\rightarrow G, \tr)$, and a framed Reidemeister pair $\Phi^x=(\phi^x, \psi^x)$,  where $\phi^x, \, \psi^x: G \times G \rightarrow E$, such that:
\begin{eqnarray*}
\d(\phi^x(L,M)) &= & [Mx^{-1},Lx^{-1} ],\\
\d(\psi^x(L,M)) &= &  [xML^{-1}x^{-1}Lx^{-1},Lx^{-1}]^{-1}.
\end{eqnarray*}
Recall that $[xML^{-1}x^{-1}Lx^{-1},Lx^{-1}]^{-1}=[L,M][ML^{-1},x]$, {in any group $G$.}
\end{Definition}

Let $(E \ra{\de} G \to \{1\},\{,\})$ be a braided crossed module. Following the notation of Definition \ref{2cmlg}, let $e \tr' k=k\{\de(k)^{-1},e\}$, for $e \in G$ and $k \in E$, thus $(E \ra{\de} G, \tr')$ is a crossed module.
\begin{Theorem}\label{el2}
 Given $x \in E$ define $\Phi^x=(\phi^x,\psi^x)$ as:
\begin{eqnarray*}
\phi^x(L,M) &= & \{Mx^{-1},Lx^{-1} \},\\
\psi^x(L,M) &= & \{xML^{-1}x^{-1}Lx^{-1},Lx^{-1}\}^{-1}.
\end{eqnarray*}
Then $\Phi^x$ is a framed Reidemeister pair{, and constitutes a} framed Eisermann lifting for $(G,x)$. 
\end{Theorem}
\begin{Proof}
 The second assertion  trivially follows from \eqref{pcbcm}. As for the first, the invariance under the Reidemeister 3 move (Definition \ref{rp}) is equivalent to equation \eqref{r3sharp}. The invariance under Reidemeister 2 is immediate. To show (i) of Definition \eqref{frp}: given $Z\in G$, the equation $\de (\phi^x(A,Z)A=Z$ implies $[Zx^{-1}, Ax^{-1}]A=Z$, giving the unique solution $A=f(Z)=Z$. Furthermore (ii) of Definition \eqref{frp} holds since $g(A)=\de(\psi^x(A,A)^{-1})A=A$.
\end{Proof}

\subsection{Homotopy interpretation of the liftings of the Eisermann Invariant}
To give a homotopy interpretation of the liftings of the Eisermann invariant, Theorems \ref{el1} and \ref{el2}, and Corollary \ref{liftextension2}, we  need to introduce the notion of non-abelian tensor product and wedge product of groups, {due to Brown and Loday {\cite{BrL0,BrL}}}. 

\subsubsection{The non-abelian tensor product and wedge product of groups}
Let $G$ be a group. Following Brown and Loday {\cite{BrL0,BrL}} and Brown, Johnson and Robertson \cite{BJR}, we define the tensor product $G\otimes G$ (a special case of the tensor product of two groups $G\otimes H$) as the group generated by the symbols $g\otimes h$, {where $g,h\in G$,} subject to the relations, $\forall g,h,k \in G$:
\begin{eqnarray}
 gh\otimes k & = &(ghg^{-1} \otimes gkg^{-1})\,\, (g\otimes k), \label{tp1}\\
 g\otimes hk & = & (g\otimes h)\,\, (hgh^{-1} \otimes hkh^{-1}). \label{tp2}
\end{eqnarray}
{It was shown by Brown and Loday in \cite{BrL}  that if $G$ is finite, then so is $G\otimes G$. That $G \otimes H$ is finite if $G$ and $H$ are finite was proven by Ellis in \cite{El}.} The following relations which hold in $G\otimes G$ are taken from {\cite{BrL0,BJR}}:
\begin{align}
 g^{-1}\otimes ghg^{-1} & = (g\otimes h)^{-1} = hgh^{-1} \otimes h^{-1} ,\label{tp3}\\
 (g\otimes h) \,\, (k\otimes m)\,\,  (g\otimes h)^{-1} & =  ([g,h]k[g,h]^{-1}) \otimes  ([g,h]m [g,h]^{-1}) ,\label{tp4} \\
 [g,h]\otimes k  & =  (g\otimes h)\,\, (kgk^{-1} \otimes khk^{-1})^{-1} ,\label{tp5} \\
 k \otimes [g,h] & =  (kgk^{-1} \otimes khk^{-1})\,\, (g\otimes h)^{-1} , \label{tp6} \\
[g \otimes h, k \otimes m] & =  [g,h] \otimes [k,m].  \label{tp7}
\end{align}

{A key fact about the {non-abelian} tensor product of groups} is that there is a homomorphism of groups $\de: G\otimes G \rightarrow G'=[G,G]$, defined on generators by $g\otimes h \mapsto [g,h]$, which is clearly surjective. Surjectivity also holds if we replace $G\otimes G$ by the group $G\wedge G$, obtained from $G\otimes G$ by imposing the {additional} relations: 
$$g\otimes g=1, \, \forall g\in G.$$ 
We denote the image of $g\otimes h$ in $G\wedge G$ by $g\wedge h$. Finally \cite{BJR}, there is a left action $\bullet$ by automorphisms of $G$ on $G\otimes G$ and $G\wedge G$, given by:
$$
{g\bullet (h\otimes k) = (ghg^{-1}) \otimes (gkg^{-1}) \textrm{ and } g\bullet (h\wedge k) = (ghg^{-1}) \wedge (gkg^{-1}),}
$$
and we have two crossed modules of groups: $(\de:G\otimes G\rightarrow G', \bullet)$ and  $(\de:G\wedge G\rightarrow G', \bullet)$. This fact (which is not immediate) is Proposition 2.5 of \cite{BrL}. {It also holds that (Proposition 2.3 of \cite{BrL}):}
$$\de(k) \tn e=k \,\, (e \bullet {k^{-1}}) \textrm{ and } e \tn \de(k)=(e \bullet k)\,\,k^{-1}.$$
These facts will be crucial in the proof of the following {lemma:}
\begin{Lemma} Let $G$ be a group.
 One has a 2-crossed module $G \otimes G \ra{\de} G \to \{1\}$, where the Peiffer lifting is:
$$\{n,m\}=(nmn^{-1} )\otimes n^{-1}= (m \tn n)^{-1}.$$
{The same statement holds for $G \wedge G \ra{\de} G \to \{1\}$, with}
${\{n,m\}=(m \wedge n)^{-1}=n \wedge m.}$
\end{Lemma}
\begin{Proof}
Conditions (i) and (ii) are trivial.

Note that, in the notation of Definition \ref{2cmlg}, {for each $f,g,h \in G$:}
\begin{align*} {f}\tr' \{h,g\}&=\{h,g\} \{\de(\{h,g\}^{-1}),{f}\} =\{h,g\} \{[h,g]^{-1},{f}\}=\{h,g\} \{[g,h],{f}\}\\
                           &=(g \tn h)^{-1}\,\, ({f} \tn [g,h])^{-1}=(g \tn h)^{-1}\,\,\big(({f}g{f}^{-1}) \otimes ({f}h{f}^{-1})\,\, (g\otimes h)^{-1}\big)^{-1}\\
             &=\big(({f}g{f}^{-1}) \otimes ({f}h{f}^{-1}) \big)^{-1}=\{{f}h{f}^{-1},{f}g{f}^{-1}\}.
\end{align*}
In particular {$f \bullet (g \tn h)=f \tr' (g \tn h)$ for each $f,g,h \in G$. This generalises, namely, for each $e \in G$ and $k \in G \otimes G$, we have:}
$$ e \tr' k=k \,\, \{\de(k)^{-1},e\}=k\,\, (e \tn \de(k)^{-1})^{-1}=k\,\, \big((e \bullet k^{-1})\,\,k\big)^{-1}=e \bullet k.$$
Also,  for each $e \in G$ and $k \in G \otimes G$:
$$\{e,\de(k)\}=(e\de(k)e^{-1}) \tn e^{-1}=\de(e\tr' k) \tn e^{-1} =(e \tr' k)\,\,k^{-1}.$$

{To prove (iii), given $k,l \in G \tn G$, we have, since $(\de\colon G \tn G \to G,\tr'=\bullet)$ is a crossed module:} 
$${\{\de(k),\de(l)\}=(\de(l)\tn \de(k))^{-1} =(l\,\,(\de(k) \tr l)^{-1})^{-1} =[l,k]^{-1}=[k,l].} $$
To prove (iv) note that if $e\in G$ and $k \in G \tn G$ then:
$$\{\de(k),e\}\{e,\de(k)\}=\big (e \tn \de(k)\big)^{-1} \,\, (\de(k) \tn e)^{-1}=k\,\, e \tr' k^{-1}\,\, e \tr' k\,\,k^{-1}=1.$$

Now, if $e,f,g \in G$:
\begin{align*}
 \{ef,g\}&=(efgf^{-1}e^{-1})\tn (f^{-1}e^{-1})=(efgf^{-1}e^{-1})\tn (e^{-1}\,\, \big(e f^{-1}e^{-1})\big)
\\&=\big((efgf^{-1}e^{-1})\tn e^{-1}\big) \,\,\big( (fgf^{-1})\tn  f^{-1}\big)\\
          &=\{e,fgf^{-1}\}\,\,\{f,g\}.
\end{align*}
Also:
\begin{align*}
 \{k,gh\}=\big((gh)\tn k\big)^{-1}&=\big (  (ghg^{-1} \otimes gkg^{-1})\,\, (g\otimes k)  \big)^{-1}\\
          &= (g\otimes k) ^{-1}\,\, g \tr' (h\tn k)^{-1}=\{k,g\}\,\,g \tr' \{k,h\}.
\end{align*}
This means (v) and (vi)  of Definition \ref{2cmlg} hold.

{Finally the second assertion follows trivially from the first.}
\end{Proof}

We can see, by taking inverses of \eqref{tp1} and \eqref{tp2}, that we could define the non-abelian tensor product of groups as being generated by $\{e,g\} =(g\otimes e)^{-1}$, with relations:
 $$\{ef,g\}=\{efe^{-1},ege^{-1}\}\,\, \{e,g\} \textrm{ and } \{e,fg\}=\{e,f\}\,\,\{fef^{-1},fgf^{-1}\}  .$$
These hold in any braided crossed module, by equations \eqref{ef,g}, the last equation of the definition of a 2-crossed module,  and \eqref{actform}.

The following theorem  is fully proved in \cite{BrL}, the first equation being also proved in \cite{M}. {Group homology} is taken with coefficients in $\mathbb{Z}$.
\begin{Theorem}[Brown-Loday / Miller]\label{blm}
Let $G$ be a group. One has  exact sequences:
$$\{0\} \to {H}_2(G) \to G \wedge G \ra{\de} G' \to \{1\}, $$
and
$$ \{0\} \to J_2(G) \to G \tn G \ra{\de} G' \to \{1\},$$ 
where $J_2(G)\doteq \ker \de$ can be embedded into the exact sequence:
$${H}_3(G) \to \Gamma(G^{\rm ab}) \ra{\zeta} J_2(G) \to {H}_2(G) \to \{0\}. $$
Here $\Gamma$ is the Whitehead universal quadratic  functor \cite{W2}, appearing in Whitehead's Certain Exact Sequence. As usual $G'=[G,G]$. 
\end{Theorem}
The explanation of the third exact sequence will be given in \ref{fint}.
\subsubsection{Homotopy interpretation of the unframed Eisermann Lifting}

Consider  a 2-crossed module of the form
$$
( E \ra{\de} G \ra{\d} 1, \{ \, , \, \})
$$
such that  $\de$ is surjective; {in particular if $g \in G$ then $\{g,g\}=\{\de(k),\de(k)\}=[k,k]=1_E$.}

 Let $K$ be a knot. Then it is well known, and {follows from the asphericity of the knot complement $C_K$ \cite{Pa}, combined with the fact that knot complements are homology circles \cite{BZ},} that ${H}_2(\pi_1(C_K))=\{0\}$. Therefore $$\pi_1(C_K)\wedge \pi_1(C_K)\cong [\pi_1(C_K),\pi_1(C_K)],$$
canonically.  
Let now $f\colon \pi_1(C_K) \to G$ be a group morphism. 
Define $$\hat{f}=\pi_1(C_K) \wedge \pi_1(C_K) \to E,$$ as acting on the generators $x \wedge  y$ of $\pi_1(C_K)\wedge \pi_1(C_K)\cong [\pi_1(C_K),\pi_1(C_K)]$ by:
$$\hat{f}(x \wedge y)=\{f(y),f(x)\}^{-1}. $$
Then  $\hat{f}$ respects the defining relations for the non-abelian  wedge product.  We therefore have a 2-crossed module map:
$$\big(\pi_1(C_K) \wedge \pi_1(C_K) \ra{\de} \pi_1(C_K)\to \{1\} \big) \ra{(\hat{f},f,\id)} \big(E  \ra{\de} G \to \{1\}\big).$$

Going back  to the  knot $K$, choose a base point $p \in K$. let $m_p \in \pi_1(C_K)$  and $l_p \in \pi_1(C_K)$ be the associated meridian and longitude. Then $$l_p \in [\pi_1(C_K) , \pi_1(C_K)] 	\cong \pi_1(C_K) \wedge \pi_1(C_K) .$$

Given an element $ x \in G$, we thus have a knot invariant of the form:
$$ \sum_{f \colon \pi_1(C_K) \to G \textrm{ with } f(m_p)=x} \hat{f}(l_p)\in \Z[E]. $$

\begin{Theorem}
 {Given $x\in G$, let $\Phi^x$ be the unframed Reidemeister pair derived from the braided crossed module $( E \ra{\de} G \ra{\d} 1, \{ \, , \, \})$, where $\de$ is surjective; see \ref{tuc}. 
 Let $K$ be a knot, with a base point $p$. Let $L_K$ be the associated string knot. Then:}
$$\sum_{a \in G}\langle 1_G \left | I_{\Phi^x}(L_K) \right | a \rangle = \sum_{f \colon \pi_1(C_K) \to G \textrm{ with } f(m_p)=x} \hat{f}(l_p). $$
\end{Theorem}
\begin{Proof}
 The proof is  analogous to the one for the unlifted case. We note Remark \ref{longformulanew}.
\end{Proof}

\subsubsection{Homotopy interpretation of the framed Eisermann Lifting}\label{fint}
The framed knot invariant derived from a braided crossed module {$(E\ra{\de} G \to \{1\}, \{,\})$} and $x \in G$, {by using the framed Reidemeister pair $\Phi^x$ of \ref{fel},} can also be given an explicit homotopy definition. 

Let $A$ be an abelian group. Recall that $\Gamma(A)$ is defined as being the group generated by all monoids $\gamma(a)=(a,a)$, with $a\in A$, with relations:
\begin{align*}
(-a,-a)&=(a,a)\\ (a+b+c,a+b+c)+(a,a)+(b,b)+(c,c)&=(a+b,a+b)+(b+c,b+c)+(c+a,c+a). 
\end{align*}
The group $\Gamma(A)$ is always abelian. It holds that $\Gamma(\Z)=\Z$, an isomorphism being, on generators: $$(a,a) \mapsto a^2.$$
Given a group $G$, there exists a group map: $$\zeta:\Gamma(G^{\rm ab}) \to G \otimes G,$$ given by $$g^{ab}\mapsto g \tn g.$$
Recalling Theorem \ref{blm}, if {${H}_2(G)$ and ${H}_3(G)$ are trivial, then $\zeta$ defines an isomorphism:} $$\zeta\colon \Gamma(G^{\rm ab}) \to \ker(\de\colon G \tn G \to G').$$

 {Let $K$ be a knot. Then we know that ${H}_2(\pi_1(C_K))=\{0\}={H}_3(\pi_1(C_K))$. Therefore} $$J_2(  \pi_1(C_K) )\doteq \ker\big(\de\colon  {\pi_1(C_K)} \otimes \pi_1(C_K) \to \pi_1(C_K)\big)\cong\Gamma(\pi_1(C_K)^{\rm ab})\cong \Gamma(\Z)\cong\Z.$$
Let $m$ be a meridian of $K$. Then $m^{\rm ab}$ is a generator of $\pi_1(C_K)^{\rm ab}=\Z$ and the previous identification $J_2(  \pi_1(C_K) )\cong \Z$ is made canonical through: $$a \mapsto (m\otimes m)^a, \textrm{ where } a \in \Z.$$

Let now $f\colon \pi_1(C_K) \to G$ be a group morphism. 
Define $$\hat{f}=\pi_1(C_K) \otimes {\pi_1(C_K)} \to E$$ as (in the generators $x \otimes y$ of $\pi_1(C_K)\otimes \pi_1(C_K)$):
$$\hat{f}(x \otimes  y)=\{f(y),f(x)\}^{-1}. $$
Then clearly $\hat{f}$ respects the defining relations for the non-abelian tensor product.  We therefore have a 2-crossed module map:
$$\big(\pi_1(C_K) \otimes \pi_1(C_K) \ra{\de} \pi_1(C_K)\to \{1\} \big) \ra{(\hat{f},f,\id)} \big(E  \ra{\de} G \to \{1\}\big).$$

{Going back  to the  knot $K$, considering a diagram $D$ of $K$, let $m\in \pi_1(C_K)$ be a meridian of $K$. Let $l\in  {[\pi_1(C_K),\pi_1(C_K)] }$ be the associated longitude. There are several elements of $ \pi_1(C_K) \tn  \pi_1(C_K) $ mapping to $l$ under $\de$. However the knot diagram $D$ expresses the longitude $l$ as a product of {(conjugates of)} commutators in $ \pi_1(C_K) $, by using Remark \ref{longformulanew}. The longitude $l$ can therefore be unambiguously lifted to $l_D\in \pi_1(C_K) \otimes  \pi_1(C_K) $, by switching commutators to Peiffer commutators, {and doing the usual categorical group evaluations of \eqref{ev1} and \eqref{ev2}.} This element $l_D$ is invariant under Reidemeister moves {1', 2 and 3}, by the calculation in \ref{fel}. Analogously to the unframed Eisermann lifting, we have:}
\begin{Theorem}
 Let $K$ be a knot, and $l$ a longitude of $K$. Let $L_K$ be the associated string knot. Given {$x\in G$} then, given any {diagram $D$} of $K$:
$$\sum_{a \in G}\langle 1_G \left | I_{\Phi^x}(L_K) \right | a \rangle = \sum_{f \colon \pi_1(C_K) \to G \textrm{ with } f(m)=x} \hat{f}(l_D). $$
\end{Theorem}

\end{document}